
\documentstyle{amsppt}
\input epsf

\topmatter
\baselineskip=1.1\baselineskip  
\loadmsbm
\UseAMSsymbols
\hoffset=0.75truein
\voffset=0.5truein

\def\x{\relax}
\def\reft{\relax}  
\def\refz{\relax}
\def\bul{\noindent$\bullet$\quad }

\def\cir{\noindent$\circ$\quad }
\def\heads#1{\rightheadtext{#1}}
\def\doct{\delta_{oct}}
\def\pt{\hbox{\it pt}}
\def\Vol{\hbox{vol}}
\def\sol{\operatorname{sol}}
\def\quo{\operatorname{quo}}
\def\anc{\operatorname{anc}}
\def\cro{\operatorname{crown}}
\def\vor{\operatorname{vor}}
\def\octavor{\operatorname{octavor}}
\def\dih{\operatorname{dih}}
\def\arc{\operatorname{arc}}
\def\rad{\operatorname{rad}}
\def\A{{\bold A}}
\def\squander{(4\pi\zeta-8)\,\pt}
\def\score{8\,\pt}
\def\Sfour{\Cal{\bold S}_4^+}
\def\Sminus{\Cal{\bold S}_3^-}
\def\Splus{\Cal{\bold S}_3^+}
\def\maxpi{\pi_{\max}}
\def\xiG{\xi_\Gamma}
\def\xiV{\xi_V}
\def\xik{\xi_\kappa}
\def\xikG{\xi_{\kappa,\Gamma}}

\font\twrm=cmr8
\def\DLP{\operatorname{D}_{\hbox{\twrm LP}}}
\def\ZLP{\operatorname{Z}_{\hbox{\twrm LP}}}

\def\geom{{\operatorname{g}}}
\def\anal{{\operatorname{an}}}

\def\E{{\Cal E}}

\def\was#1{\relax}

\def\diag|#1|#2|{\vbox to #1in {\vskip.3in\centerline{\tt Diagram #2}\vss} }
\def\v{\hskip -3.5pt }
\def\gram|#1|#2|#3|{
        {
        \smallskip
        \hbox to \hsize
        {\hfill
        \vrule \vbox{ \hrule \vskip 6pt \centerline{\it Diagram #2}
         \vskip #1in %
             \includegraphics{#3}\hrule }
        \v\vrule\hfill
        }
\smallskip}}

\def\samgram|#1|#2|#3|{
        {
        \smallskip
        \hbox to \hsize
        {\hfill
        \vrule \vbox{ \hrule \vskip 6pt \centerline{\it Diagram #2}
         \vskip #1in %
            \centerline{\epsfbox{#3}}
         \vskip #1in %
    \hrule }
        \v\vrule\hfill
        }
\smallskip}}


\def\today{\ifcase\month\or
    January\or February\or March\or April\or May\or June\or
    July\or August\or September\or October\or November\or December\fi
    \space\number\day, \number\year}

\title Sphere Packings IV
\endtitle

\author Thomas C. Hales\endauthor
\endtopmatter

\document

\footnote""{\hfill version -- 7/31/98, revised 1/21/02}
\footnote""{Research partially funded by the NSF.}

\bigskip
\head \refz 1. Introduction and Review\endhead

\subhead \refz 1.1. The steps\endsubhead
The Kepler conjecture asserts that no packing of spheres in
three dimensions has density greater than
$\pi/\sqrt{18} \approx 0.74048$.  This paper is one of a series
of papers devoted to the Kepler conjecture.  This series began
with \cite{I}, which proposed a line of research to prove the conjecture,
and broke the conjecture into smaller conjectural steps
which imply the Kepler conjecture.  The
steps were intended to be equal in difficulty,
 although some have emerged
as more difficult than others.
This paper completes part of the fourth step.  The main result
is Theorem 4.4.

As a continuation of \cite{F} and \cite{III}, this paper assumes
considerable familiarity with the constructions, terminology, and notation
from these earlier papers.  See \cite{F} for the definitions
of quasi-regular tetrahedra, upright and flat quarters and their
diagonals, anchors, Rogers simplices, standard clusters, standard
regions, the $Q$-system, $V$-cells, local $V$-cells,
and decomposition stars.

We will use a number of constants and functions from \cite{I} and
\cite{F}: $\doct=(\pi-2\zeta^{-1})/\sqrt8$,
$\pt=-\pi/3+2\zeta^{-1}$,
$\zeta^{-1}=2\arctan(\sqrt2/5)$, $t_0=1.255$,
$\phi(h,t)$, $\phi_0=\phi(t_0,t_0)$,
$\Gamma(S)$ is the compression, $\vor(S)$ is the
analytic Voronoi function, $\vor(\cdot,t)$ and $\vor_0=\vor(\cdot,t_0)$
are the truncated Voronoi functions, $\sigma(S)$ is the score,
$\tau(S)$ measures what is squandered by a simplex,
$\dih(S)$ is the dihedral angle along
the first edge of a simplex, $\sol$ is the solid angle, $R(a,b,c)$ is a
Rogers simplex with parameters $a\le b\le c$, and $S(y_1,\ldots,y_6)$ is a
simplex with edge lengths $y_i$ with the standard conventions from \cite{I}
on the ordering of edges.  The definitions of $\sigma(S)$ and
$\tau(S)$ are particularly involved.  The definition depends on the
structure and {\it context\/} of $S$
\cite{F.3}.

At the heart of this approach is a geometric structure, called
the {\it decomposition star}, constructed around the center of each
sphere in the packing.  A function $\sigma$,
called the {\it score}, is defined
on the space of all decomposition stars.  An upper
 bound of $\score\approx 0.4429$ on the
score implies the Kepler conjecture \cite{F, Proposition 3.14}.
A second function $\tau$,
measuring what is {\it squandered}, is defined on the space of
decomposition stars.  If a decomposition
star squanders
more than $\squander\approx 14.8\,\pt\approx 0.819$,
then it scores less than $\score$.

\subhead 1.2 Exceptional regions\endsubhead

A standard region is said to be {\it exceptional\/} if it is
not a triangle or quadrilateral.
The vertices of the packing of height at most $2.51$ that
are contained in the closed cone over the standard region
are called its {\it corners}.

The results of this paper are based on a number of inequalities listed
in the appendix.  These inequalities are grouped into collections
denoted $\A_i$.

\subhead 1.3 Organization of this paper\endsubhead

Lemmas, Remarks, Propositions, and so forth, are numbered
according to the following conventions.  The paper is divided into
five sections and Appendices.  Each section is divided into a
number of subsections. The reference $n.m$ refers to Subsection
$m$ of Section $n$, or more briefly, Section $n.m$.  For instance,
this is Section $1.3$. In general, lemmas, remarks, and so forth,
are numbered according to the subsection in which they appear.
Thus, Proposition $3.7$ is the unique Proposition in Section
$3.7$. When more than one Lemma appears in a subsection, they are
numbered consecutively.  Thus, the three lemmas of Section $3.8$
are Lemmas $3.8.1$, $3.8.2$, and $3.8.3$.

Appendix $1$ contains long listings of inequalities that are used
throughout the paper.  These inequalities are grouped into $24$ sections
$\A_1,\ldots,\A_{24}$.  Each inequality is labeled with an integer
$1,\ldots,k$ and a unique nine-digit identification number.  In
the body of the paper, each inequality is identified by its Section
and integer label.  Thus, Inequality $\A_8.4$ is the fourth inequality
in Section $\A_8$ of Appendix $1$.  The nine-digit identification code
is used to identify the inequality in the archive of computer code
that was used to test and prove the inequality.  These numeric codes
make it easy to locate computer files that deal with a particular
inequality.

\bigskip
\head \refz 2. The fine decomposition\endhead

\subhead\refz 2.1. Overview of the Fine Decomposition\endsubhead

In Section 2, we define a decomposition of each local $V$-cell
$V_P$, called the {\it fine decomposition\/} of $V_P$.  Let $V$ be
the $V$-cell at the origin.  Let $P$ be an exceptional
cluster.  Recall that the part of $V$ in the cone
$C(P)$ over $P$ is called the local
$V$-cell $V_P$.  Let $V_P(t_0)$ be the intersection of $V_P$ with
the ball $B(t_0)$ of radius $t_0 = 1.255$.  We write $V_P$ as
the disjoint union of $V_P(t_0)$ and its complement $\delta_P$.

Let $v$ be an enclosed vertex of height between 2.51 and $2\sqrt{2}$.
Assume that there is an upright quarter in the $Q$-system with
diagonal $(0,v)$.
We call $(0,v)$ an {\it upright diagonal}.
  We will define $\delta_P(v)\subset\delta_P$.
It will be a subset of a set of the form $C(D_v)\cap\delta_P$ for
some subset $D_v$ of the unit sphere.  The sets $D_v$ will be defined
so as not to overlap one another for distinct $v$.  Then the sets
$\delta_P(v)$ do not overlap one another either.   We will give an
explicit formula for the volume of $\delta_P(v)$.

We will define a set $\Cal S$ of simplices in $C(V_P)$.  The vertices
of the simplices will be vertices of the packing, and their edges
will have length at most $2\sqrt{2}$.  The sets $C(S)$,
for distinct  $S\in\Cal S$,
 will not overlap.
Over a simplex $S\in\Cal S$,
  the local $V$-cell will be
truncated at a radius $t_S\ge t_0$.
  After defining the constants $t_S$, we will set
$$V_S(t_S) = C(S)\cap B(t_S)\cap V_P.$$
If $V_P\cap C(S)\subset B(t_S)\subset B(t'_S)$, then $V_S(t_S)=V_S(t'_S)$.

Since $t_S\ge t_0$, the sets $V_S(t_S)$ and $\delta_P$
may overlap.
Nevertheless, we will show that $V_S(t_S)$ does not overlap any
$\delta_P(v)$.  Let $\tilde V_P(t_0)$ be the set of
points in $V_P(t_0)$ that do not lie in $C(S)$,
$S\in\Cal S$.  We will derive an explicit formula for the volume
and score of $\tilde V_P(t_0)$.

In $V_P$, there are nonoverlapping sets
$$\delta_P(v),\quad   V_S(t_S),\quad \tilde V_P(t_0).$$
Let $\delta_P'$ be the complement in $V_P$ of the union of these sets.
These sets give a decomposition of $V_P$,
called the fine decomposition of the local $V$-cell $V_P$.
Corresponding to the fine decomposition is a formula for
the score $V_P$ of the form
$$\sigma(V_P) =
    \sigma(\tilde V_P(t_0))
    + \sum_{\Cal S}\sigma(V_S(t_S))
    -\sum_v 4\doct\Vol(\delta_P(v))-4\doct\Vol(\delta_P').$$
Since $\Vol(\delta_P')\ge0$, we obtain an upper bound on the score
of $V_P$ by dropping the rightmost term.

\subhead\refz 2.2. $V$-cells\endsubhead

Let $\Cal Q_0$ be the set of simplices
in the $Q$-system
with a vertex at the origin.
If $x$ lies in the Voronoi cell at the origin, but not in the $V$-cell,
then either $x$ belongs to a simplex in $\Cal Q_0$ or $x$ belongs
to the protruding tip from a simplex in $\Cal Q_0$.  In either case,
$x\in C(S)$, for some $S\in \Cal Q_0$.
Consequently, the part of the Voronoi cell over the complement
of $C(S)$, for all $S\in \Cal Q_0$,
is contained in the $V$-cell.

\subhead\refz 2.3. The set $\delta_P(v)$\endsubhead

Let $(0,v)$ be the diagonal of an upright quarter in $\Cal Q_0$ and in
    the cone over $P$.
We define $\delta_P(v)\subset C(D_v)\cap \delta_P$ for an appropriate
subset $D_v$ of the unit sphere.

Let $D_0$ be the spherical cap on the unit sphere, centered along
    $(0,v)$ and having arcradius $\theta$, where
    $\cos\theta = |v|/(2\eta_0(|v|/2))$,
    and $\eta_0(h)=\eta(2h,2,2.51)$.  The area of $D_0$ is
    $2\pi(1-\cos\theta)$.
Let $v_1,\ldots,v_k$ be the anchors around $(0,v)$ indexed cyclically.
    The projections of the edges $(v,v_i)$ (extended as necessary)
    slice the spherical cap into $k$ wedges $W_i$, between
    $(v,v_i)$ and $(v,v_j)$, where $j\equiv i+1\mod k$, so that
    $D_0 =\cup W_i$.

Let $\Cal W$ be the set of wedges $W=W_i$  such
that either

\noindent(1) $W$ occupies more than half the spherical cap (so that
    its area is at least $\pi(1-\cos\theta)$), or

\noindent(2) $|v_i-v_j|\ge 2.77$,
    $\rad(0,v,v_i,v_j)> \eta_0(|v|/2)$, and the
    circumradius of $(0,v_i,v_j)$ or $(v,v_i,v_j)$ is
    $\ge\sqrt2$.

Fix $i,j$, with $j\equiv i+1\mod k$.
If $W = W_i$ is a wedge in $\Cal W$, let $(0,v_i,v)^\perp$ be
the plane through the origin and the circumcenter of $(0,v_i,v)$, perpendicular
to $(0,v_i,v)$.  Skip the following step if the circumradius
of $(0,v_i,v)$ is greater than $\eta_0(|v|/2)$, but if the circumradius
is at most this bound, let $c_i$ be the intersection of $(0,v_i,v)^\perp$
with the circular boundary of $W$.  Extend $W$ by adding to $W$ the
spherical triangle with vertices the projections of
$v$, $v_i$, and $c_i$.  Similarly, extend $W$ with the
triangle from $(v,v_j,c_j)$, if the circumradius of $(0,v_j,v)$ permits.
(An example of this is illustrated in F.4.6.)  Let $W^e$ be extension
of the wedge obtained by adding these two spherical triangles.

We will define $\delta_P(W^e)\subset C(W^e)\cap \delta_P$.  Then
$\delta_P(v)$ is defined as the union of $\delta_P(W^e)$, for $W\in\Cal W$.
Let
$$E_w = \{x : 2 x\cdot w \le w\cdot w\},$$
for $w = v,v_i,v_j$.  These are half-spaces bounding the Voronoi cell.
Set $E_\ell = E_{v_\ell}$.

If (2) holds, we let $c$ be the projection of the circumradius of
$(0,v_i,v_j,v)$ to the unit sphere.  The arclength from
$c$ to the projection of $v$ is $\theta'$, where
$$\cos\theta' = |v|/(2\rad)<|v|/(2\eta_0) = \cos\theta.$$
We conclude that $\theta'>\theta$ and $c$
does not lie in $D_0$.

In both cases (1) and (2) set
$$\delta_P(W^e) = (E_v\cap E_i\cap E_j \cap C(W^e))\setminus B(t_0).$$
Observe that
$$E_v\cap E_i\cap E_j \cap C(W^e)$$
is the union of four Rogers simplices
$$R(|w|/2,\eta(0,v,v_\ell),\eta_0(|v|/2)),\quad w = v,v_\ell,\quad
    \ell = i,j$$
and a conic wedge over $W$ between $c_i$ and $c_j$.
(The inequality $\theta'>\theta$ implies that
the Rogers simplices do not overlap.)

\proclaim {Lemma}  $\delta_P(W^e)\subset V_P$.
\endproclaim

\demo{Proof} First assume for a contradiction that
some part of the wedge between $c_i$ and $c_j$ overlaps
the $V$-cell at some other vertex $v'$.  This forces
$\eta(0,v,v')<\eta_0(|v|/2)$, and $v'$ must then be an anchor.
But for anchors, the separation of $V$-cells
has been achieved
by the half-spaces $E_*$.

Now suppose one of the Rogers simplices along $(0,v,v_i)$ overlaps
the $V$-cell at some vertex $v'$.  If $v'$ lies on the opposite
side of the plane $(0,v,v_i)$ from $W^e$, then in order to meet, the
face $(0,v,v_i)$ of $(0,v,v_i,v')$ must have negative orientation.
This forces
$(0,v,v_i,v')$ to be a quarter \cite{F}.
It is in the $Q$-system because one
and hence all quarters along $(0,v)$ lie in the $Q$-system.  Thus
any protruding tip from $v'$ is reapportioned among neighboring $V$-cells,
so that such a point of $\delta_P(W^e)$ lies in the $V$-cell.

Take $v'$ and $W^e$ to lie on the same side of $(0,v,v_i)$.
To overlap, the circumradius of $(0,v,v_i,v)$ must be less than $\eta_0(|v|/2)$.
Then $\eta(0,v,v')<\eta_0(|v|/2)$, forcing $v'$ to be an anchor.
Since $v_i$ and $v_j$ are chosen to be consecutive, we find that
$v'=v_j$.  But then the condition (2) gives the contradiction
$\rad(0,v,v_i,v')\ge\eta_0(|v|/2)$.
\qed
\enddemo

\subhead\refz 2.4. Overlap\endsubhead

\proclaim{Proposition}  The sets $\delta_P(W^e)$ do not overlap.
\endproclaim

\demo{Proof}  This is clear for two sets around the same vertex $v$.
In general, this follows from the fact that the sets $W^e$ do not
overlap on the unit sphere.  We use the faces of the $V$-cell
to separate them.  In the notation of Sections \refz 2.1-2.3,
the part of the wedge $W$ between $c_i$ and $c_j$ lies under the
face of the $V$-cell associated with $v$, the vertex used to
construct $W$.  Hence, these pieces do not overlap at different vertices.
Similarly, two of the Rogers simplices lie under the face of the
$V$-cell associated with $v$.  The remaining two Rogers simplices
lie under the faces of $V$-cells of two of the anchors of $v$.
A vertex $v_i$ may be the anchor of more than one upright diagonal $(0,v)$
and $(0,v')$.
Nevertheless, the corresponding Rogers simplices do not overlap because each
Rogers simplex for $W^e$ at $v$ will lie under the triangular part of the face
determined by $v_i/2$ and the edge of the $V$-face dual to the
triangle $(0,v_i,v)$, and the Rogers simplex for $W^{\prime e}$ at $v'$
will lie under a corresponding triangular part of the face.  These
triangles do not overlap, so the extended wedges cannot either.
\qed
\enddemo

Suppose that the faces of the $V$-cell
dual to two vertices $v_1$ and $v_2$
(of height at most $2\sqrt{2}$) share
an edge.  On the face dual to $v_1$, we take the triangle formed
by $v_1/2$ and the common edge, and call it the $(v_1,v_2)$-triangle.
(Since $|v_1|\le2\sqrt{2}$, $v_1/2$ lies on the face dual to $v_1$.)
The proof shows that the set $\delta_P(v)$ lies under the face dual to
$v$ or under the $(v_i,v)$-triangles of anchors $v_i$ of $v$.

\subhead\refz 2.5. Some simplices\endsubhead

We consider three types of simplices $S_A$, $S_B$, $S_C$.  Each type
has its vertices at vertices of the packing.  The
edge lengths of these simplices are at most $2\sqrt{2}$.

$S_A$.  This family consists of simplices $S(y_1,\ldots,y_6)$ whose
    edge lengths satisfy
    $$y_1,y_2,y_3\in[2,2.51],\quad y_4,y_5\in[2.51,2.77],
    \quad y_6\in[2,2.51],\quad \text{and }\eta(y_4,y_5,y_6)<\sqrt{2}.$$
(These conditions imply $y_4,y_5<2.697$, because $\eta(2.697,2.51,2)>\sqrt2$.)

$S_B$.  This family consists of certain flat quarters that are part of
    an isolated pair of flat quarters.
    It consists of those satisfying $y_2,y_3\le 2.23$,
    $y_4\in[2.51,2\sqrt{2}]$.

$S_C$.  This family consists of certain
    simplices $S(y_1,\ldots,y_6)$ with
    edge lengths satisfying
    $y_1,y_4\in[2.51,2\sqrt{2}]$, $y_2,y_3,y_5,y_6\in[2,2.51]$.
    We impose the condition that the first edge is the
    diagonal of some upright quarter in the $Q$-system,
    and that the upper endpoints
    of the second and third edges (that is, the second and third
    vertices of the simplex) are consecutive anchors of this
    diagonal.
    We also assume that $y_4\le 2.77$, or that
    both face circumradii
    of $S$ along the fourth edge are
    at most $\sqrt{2}$.

\proclaim{Lemma}  If a vertex $w$ is enclosed over a simplex $S$ of type
$S_A$, $S_B$, or $S_C$, then its height is greater than 2.77.  Also,
$(0,w)$ is not the diagonal of an upright quarter in the $Q$-system.
\endproclaim

\demo{Proof}
In case $S_A$, $\eta(y_4,y_5,y_6)<\sqrt{2}$, so an enclosed vertex
must have height greater than $2\sqrt{2}$.  It is too long to be the
diagonal of a quarter.

In case $S_B$, we use the fact that the isolated quarter does not
overlap any quarter in the $Q$-system.  We recall that a function
$\Cal E$, defined in \cite{F}, measures the distance between opposing
vertices in a pair of simplices sharing a face.
An enclosed vertex has length
at least
    $$\Cal E(S(2,2,2,2\sqrt{2},2.51,2.51),2.51,2,2)>2.77.$$
By the symmetry of isolated quarters, this means that the diagonal of
a flat quarter must also be at least $2.77$.

In case $S_C$, the same calculation gives that the enclosed vertex $w$
has height
at least $2.77$.  Let the simplex $S$ be given by
$(0,v,v_1,v_2)$, where $(0,v)$ is the upright diagonal.
By Lemma F.1.5, $v_1$ and $v_2$
are anchors of $(0,w)$.  The edge
between $w$ and its anchor cannot
cross $(v,v_i)$ by Lemma F.1.3.
(Recall that two sets are said to {\it cross\/}
if their projections overlap.)
 The distance between $w$ and
$v$ is at most $2.51$ by Lemma F.1.9.   If $(0,w)$ is the diagonal
of an upright quarter, the quarter takes the form
$(0,w,v_1,v_3)$, or $(0,w,v_2,v_3)$ for some $v_3$,
by Lemma F.1.8.  If both
of these are quarters, then the diagonal $(v_1,v_2)$ has four anchors
$v$, $w$, $0$, and $v_3$.  The selection rules for the $Q$-system
place the quarters around this diagonal in the $Q$-system.
So neither $(0,w,v_1,v_3)$ nor $(0,w,v_2,v_3)$ is in the
$Q$-system.
Suppose that $(0,w,v_1,v_3)$ is a quarter, but that
$(0,w,v_2,v_3)$ is not.  Then
$(0,w,v_1,v_3)$ forms an isolated pair with $(v_1,v_2,v,w)$.  In either case
the quarters along $(0,w)$ are not in the $Q$-system.
\qed
\enddemo

\proclaim{Remark}  The proof of
this lemma does not make use of all the
hypotheses on $S_C$.  The conclusion holds for any
simplex $S(y_1,\ldots,y_6)$,
with $y_1,y_4\in[2.51,2\sqrt{2}]$, $y_2,y_3,y_5,y_6\in[2,2.51]$.
\endproclaim

\subhead \refz 2.6. Disjointness \endsubhead

Let $S=(0,v_1,v_2,v_3)$ be a simplex of type $S_A$, $S_B$, or $S_C$.
An edge $(v_4,v_5)$
of length at most $2\sqrt{2}$ such that $|v_4|,|v_5|\le 2.51$
cannot cross two of the edges $(v_i,v_j)$ of $S$.  In fact,
it cannot cross any edge $(v_i,v_j)$ with $|v_i|,|v_j|\le 2.51$
by Lemma F.1.6.  The only possibility is that the edge $(v_4,v_5)$
crosses the two edges with endpoint $v_1$, with $|v_1|\ge2.51$
in case $S_C$.  But this too is impossible by Lemma F.1.8.

Similar arguments show that the same conclusion holds
for an edge $(v_4,v_5)$ of length at most $2.51$ such that
$|v_4|\le2.51$, $v_5\le2\sqrt{2}$.  The only additional fact that
is needed is that $(v_4,v_5)$ cannot cross the edge between
the vertex $v$ of an upright diagonal $(0,v)$ and an anchor \cite{F.1.3}.

Now take two simplices $S$, $S'$, each of  type $S_A$, $S_B$, $S_C$,
 or a quarter in the $Q$-system.

\proclaim{Lemma} $S$ and $S'$ do not overlap.
\endproclaim

\demo{Proof}
We remark that we are tacitly assuming that the standard region
is exceptional, and we exclude the case of  conflicting diagonals
in a quad cluster.
We claim that no vertex $w$ of $S$ is enclosed over $S'$.
Otherwise, $w$ must have height at least $2.51$, so that
$(0,w)$ is the diagonal of an
upright in the $Q$-system, and this is contrary to Lemma
\refz 2.5.
Similarly, no vertex of $S'$ is enclosed over $S$.

Let $(v_1,v_2)$ be an edge of $S$ crossing an edge
$(v_3,v_4)$ of $S'$.  By the preceding remarks,
neither of these edges can cross two edges of the other simplex.
The endpoints of the edges are not enclosed over the other simplex.
This means that one endpoint of each edge $(v_1,v_2)$ and $(v_3,v_4)$
is a vertex of the other simplex.  This forces $S$ and $S'$ to have three
vertices in common, say $0$, $v_2$, and $v_3$.  We have
$S=(0,v_1,v_3,v_2)$ and $S'=(0,v_3,v_2,v_4)$.
If $|v_2|\in[2.51,2\sqrt{2}]$, then we see that the anchors $v_3$, $v_4$
of $(0,v_2)$ are not consecutive.  This is impossible for simplices
of type $S_C$ and upright quarters.  Thus, $v_2$ and $v_3$ have height
at most $2.51$.  We conclude, without loss of generality,
that $|v_4|\in[2.51,2\sqrt{2}]$ and $|v_1-v_2|\ge 2.51$ \cite{F.1}.

The heights of the vertices of $S$ are at most $2.51$, so it has type
$S_A$ or $S_B$, or it is a flat quarter in the $Q$-system.
If $S'$ is an upright quarter in the $Q$-system, then it does not
overlap an isolated quarter or a flat quarter in the $Q$-system,
so $S$ has type $S_A$.
This imposes the contradictory constraints on $S_A$
$$2.77\ge |v_1-v_2|\ge\Cal E(S(2.51,2,2,2\sqrt{2},2.51,2.51),2,2,2)>2.77.$$
Thus $S'$ has type $S_C$.  This forces $S$ to have type $S_A$.  We reach
the same contradiction  $2.77\ge \Cal E>2.77$.
\qed
\enddemo

\subhead \refz 2.7 Separation of simplices of type $S_A$ \endsubhead

Let $V_S = V_P\cap C(S)$, for a simplex $S$ of type $S_A$, $S_B$, or $S_C$.
We truncate $V_S$ to $V_S(t_S)$ by intersecting $V_S$ with a
ball of radius $t_S$.  The parameters $t_S$ depend on $S$.

If $S$ has type $S_A$, we use $t_S=+\infty$ (no truncation).  If
$v$ is enclosed over $S=(0,v_1,v_2,v_3)$, then since
$\eta(v_1,v_2,v_3)<\sqrt{2}$,
the face $(v_1,v_2,v_3)$ has positive orientation for $S$ and $(v,v_1,v_2,v_3)$.
This implies that the $V$-cells at $v$ and $0$ do not intersect, and
there is no need to truncate.  If a simplex adjacent to $S$ has
negative orientation along a face shared with $S_A$, then it must
be a quarter $Q=(0,v_4,v_1,v_2)$ \cite{F.2.2} or quasi-regular tetrahedron.
It cannot be an isolated quarter because
of the edge length constraint $2.77$ on simplices of type $S_A$.
If it is in the $Q$-system, it does not interfere with the $V$-cell over $S$.
Assume that it is not in the $Q$-system.  There must be a conflicting
diagonal $(0,w)$, where $w$ is enclosed over $Q$.  ($w$ cannot be
enclosed over $S$ by results of Lemma \refz 2.6.)  This shields the
$V$-cell at $v_4$ from $C(S)$ by the two faces $(0,w,v_1)$ and $(0,w,v_2)$
of quarters in the $Q$-system.

This shows that nothing external to a simplex of type $S_A$ affects
the shape of
$V_S(t_S)$, so that $V_S(t_S)$ can be computed from $S$ alone.
Similarly, $V_S(t_S)$ does not influence the external geometry, since
all faces have positive orientation.

We also remark that $V_S(t_S)$ does not overlap any of the sets
$\delta_P(v)$.  This is evident because the two types of sets lie
under the faces of
$V$-cells associated with different vertices of the packing.
A set $\delta_P(v)$ lies under the face of the $V$-cell dual to $v$
or under the $(w,v)$-triangles of anchors $w$ of $V$.  But
$V_S(t_S)$ lies under the $(v_i,v_j)$-triangles, for the
edges $(v_i,v_j)$ of $S$.  (See Section \refz 2.4.)

Our justification that $V_S(t_S)$ can be treated as an independently
scored entity is now complete.

\subhead \refz 2.8. Separation of simplices of type $S_B$\endsubhead

If $S(y_1,\ldots,y_6)$ has type $S_B$, we label vertices so that the
diagonal is the fourth edge, with length $y_4$.
We set $t_S=1.385$.
The calculation of $\Cal E$ in Section 2.5 shows that any enclosed
vertex over $S$ has height at least $2.77=2t_S$.

    Vertices outside $C(S)$ cannot affect the shape of $V_S(t_S)$.  In
fact, such a vertex $v'$ would have to form a quarter or quasi-regular
tetrahedron with a face of $S$.  The $V$-cell at $v'$ cannot
meet $C(S)$ unless it is a quarter that is not in the $Q$-system.
But by definition, an isolated quarter is not adjacent (along
a face along the diagonal) to any other
quarters.

    To separate the scoring of $V_S(t_S)$ from the rest of the
standard cluster, we also show that the terms of Formula \refz F.3.5
for $V_S(t_S)$ lie in the cone $C(S)$.  This is more than a formality
because $S$ can have negative orientation along the face $F$ formed by
the origin and the diagonal (the fourth edge).

Let $\arc(a,b,c)=\arccos((a^2+b^2-c^2)/(2a b))$ be the angle
opposite the edge of length $c$ in a triangle with sides $a$, $b$, $c$.
    Let $\beta_\psi(y_1,y_3,y_5)\in[0,\pi/2]$ be defined by the
equations
$$\align
    \cos^2\beta_\psi &= (\cos^2\psi-\cos^2\theta)/(1-\cos^2\theta),
        \hbox{ for }\psi\le\theta,\\
    \theta &= \arc(y_1,y_3,y_5).\\
\endalign
$$
If we form a triangle $(0,v_1,v_3)$, where $|v_1|=y_1$, $|v_3|=y_3$,
$|v_1-v_3|=y_5$, then $\theta$ is the angle at the origin between
$v_1$ and $v_3$.  If we place a spherical cap of arcradius $\psi$
on the unit sphere centered along $(0,v_1)$, then the angle along $(0,v_3)$
between the plane $(0,v_1,v_3)$ and the plane tangent to the spherical
cap passing through $(0,v_3)$ is $\beta_\psi(y_1,y_3,y_5)$.

Let $S=(0,v_1,v_2,v_3)$, where $v_i$ is the endpoint of the $i$th edge.
We establish that the conic and Rogers terms of Formula F.3.5 lie
over $C(S)$ by showing that
$\beta_\psi(y_1,y_3,y_5) < \dih_3(S(y_1,\ldots,y_6))$, where $\dih_3$
is the dihedral angle along the third edge.
We use $\cos\psi = y_1/2.77$ and assume $y_2,y_3\in[2,2.23]$.
See $\A_1$.

    The reasons given in Section \refz 2.7
for the disjointness of $\delta_P(v)$ and
$V_S(t_S)$ apply to simplices of type $S_B$ as well.
This completes the justification that $V_S(t_S)$ is an object
that can be treated in separation from the rest of the local $V$-cell.

\subhead \refz 2.9. Separation of simplices of type $S_C$\endsubhead

If $S(y_1,\ldots,y_6)$ is of type $S_C$, we label vertices so that
the upright diagonal is the first edge.  We use $t_S =+\infty$
(no truncation).   Each face of $S$ has positive orientation by
F.2.2.  So $V_S(t_S)\subset S$.

    Vertices outside $S$ cannot affect the shape of $V_S(t_S)$.  Any
vertex $v'$ would have to form a quarter along a face of $S$.  If
the shared face lies along the first edge, it is a quarter $Q$ in
the $Q$-system, because one and hence all quarters along this edge
are in the $Q$-system.  If the shared face lies along the fourth
edge, then its length is at most $2.77$, so that the quarter cannot
be part of an isolated pair.  If it is not in the $Q$-system, there
must be a conflicting diagonal.  The two faces along this conflicting
diagonal of the adjacent pair in the $Q$-system (that is, the pair
taking precedence over $Q$ in the $Q$-system) shield the $V$-cell at $v'$
from $S$.

    The reasons given in Section \refz 2.7 for the disjointness of $\delta_P(v)$
and $V_S(t_S)$ apply to simplices of type $S_C$ as well.
This completes the justification that $V_S(t_S)$ is an object
that can be treated in separation from the rest of the local $V$-cell.

\subhead \refz 2.10. Simplices of type $S'_C$\endsubhead

We introduce a small variation on simplices of type $S_C$, called
type $S'_C$.  We define a simplex $(0,v,v_1,v_2)$ of type $S'_C$
to be one satisfying the following conditions.
    (1) The edge $(0,v)$ is an upright diagonal of an upright quarter
        in the $Q$-system.
    (2) $|v_2|\in[2.45,2.51]$.
    (3) $|v-v_2|\in [2.45,2.51]$.
    (4) The edge $(v_1,v_2)$
is a diagonal of a flat quarter with face $(0,v_1,v_2)$.

On simplices $S$ of type $S'_C$, we label vertices so that the
upright diagonal is the first edge.  We use $t_S=+\infty$ (no
truncation).  Each face of $S$ has positive orientation by \refz
F.2.2. So $V_S(t_S)\subset S$.

Simplices of type $S'_C$ are separated from quarters in the
$Q$-system and simplices of types $S_A$ and $S_B$ by procedures
similar to those described for type $S_C$.  The following lemma is
helpful in this regard.

\proclaim{Lemma}
 The flat quarter along the face $(0,v_1,v_2)$ is
in the $Q$-system.
\endproclaim

\demo{Proof}
$${\Cal E}(S(2,2,2.45,2\sqrt2,2t_0,2t_0),2,2,2)>2\sqrt2,$$
so nothing is enclosed over the flat quarter.
$$\E(S(2,2,2,2\sqrt2,2t_0,2t_0),2t_0,2.45,2)> 2\sqrt2,$$
so no edge between vertices of the packing can cross inside the
anchored simplex. This implies that the flat quarter does not have
a conflicting diagonal and is not part of an isolated pair.
    \qed\enddemo

Similar arguments show that there is not a simplex with negative
orientation along the  top face of $S$.

\subhead \refz 2.11. Scoring \endsubhead

The construction of the fine decomposition of the
local $V$-cell $V_P$ is now complete.
It consists of the pieces

    \bul $\delta_P(v)$, for each diagonal $(0,v)$ of an upright quarter
        in the $Q$-system,

    \bul truncations of Voronoi pieces $V_S(t_S)$ for simplices of type
        $S_A$, $S_B$, or $S_C$ ($S'_C$), over $C(P)$,

    \bul $\tilde V_P(t_0)$, the truncation at $t_0$ of all parts of
        $V_P$ that do not lie in any of the cones $C(S)$ over
        simplices
        of type $S_A$, $S_B$ or $S_C$,

    \bul $\delta_P'$, the part not lying in any of the preceding.

By the results of Sections \refz 2.7--2.9, the score of $V_P$ can be broken
    into a corresponding sum,
$$\align
\sigma(P) &= \sum_Q \sigma(Q) + \sigma(V_P),
                \hbox{ for quarters $Q$ in the $Q$-system,}\\
\sigma(U) & = 4(-\doct\Vol(U) + \sol(U)/3),\quad
    \hbox{ where }U = V_P, V_S(t_S), \tilde V_P(t_0),\\
\sigma(V_P) &= \sigma(\tilde V_P(t_0))+  \sum_{S_A,S_B,S_C} \sigma(V_S(t_S))
        - \sum_v 4\doct\Vol(\delta_P(v)) - 4\doct\Vol(\delta'_P).\\
\endalign
$$

By dropping the final term, $4\doct\Vol(\delta'_P)$, we obtain an
upper bound on $\sigma(V_P)$.  Because of the separation results
of Sections \refz 2.7--2.9, we may score $\tilde V_P(t_0)$ by the Formula
\refz F.3.7.
Bounds on the score of simplices of type $S_B$ appear in
$\A_1$.

\heads{}
\vfill\eject
\head \refz 3.  Upright Quarters\endhead

\subhead \refz 3.1. Definitions\endsubhead

Fix an exceptional cluster $R$.
Throughout this paper, we assume that $R$ lies on
a star of score at
least $\score$.
It is to be understood, when we
say that a standard region does not exist, that we mean that
there exists no such region on any star scoring more
than $\score$.

\heads{\refz 3. Upright Quarters}

In Section 3, we discuss
how to eliminate many cases of upright diagonals.
The results are summarized at the
end of the section (\refz 3.10).

If $R$ is a standard cluster or region,
we write $V_R(t)$
for the intersection of the local $V$-cell $V_R$ with a ball
$B(t)$, centered at the origin, of radius $t$.  We generally
take $t=t_0=1.255=2.51/2$.
If $(0,v)$, of length between 2.51 and $2\sqrt{2}$, is not the
diagonal of an upright quarter in the $Q$-system,
then $v$ does not affect
the truncated cell $V_R(t_0)$
and may be disregarded.  For this reason we confine our attention
to upright diagonals, which by definition lie along an upright
quarter in the $Q$-system.

\subhead \refz 3.2. Truncation\endsubhead

We say that an upright diagonal $(0,v)$
can be {\it erased with penalty $\pi_0\ge0$}, if we have, in terms of the
fine decomposition,
$$\sum_Q\sigma(Q) +\sum_S\sigma(V_S(t_S)) - 4\doct\Vol(\delta_P(v))<
    \pi_0+\sum_Q\vor_0(Q) +\sum_S\vor_0(S).$$
Here the sum over $Q$ runs over the upright quarters around $(0,v)$.
Their scores $\sigma(Q)$ are context-dependent (see \refz F.3.)
The simplices $S$ are those along $(0,v)$ of type $S_C$ in the
fine decomposition.  We define their score $\sigma(V_S(t_S))$
as in Section \refz 2.  Also, $\delta_P(v)$ is the piece
of the fine decomposition
defined in Section \refz 2.   The right-hand side is scored by
the truncation of the Voronoi function, Formula \refz F.3.7.
When we erase without mention of a penalty, $\pi_0=0$ is assumed.

If the diagonal can be erased, an upper bound on the score is obtained
by ignoring the upright diagonal and all of the structures around
it coming from the fine decomposition, and switching to the
truncation at $t_0$.   Section 3 shows that various
vertices can be erased, and this will greatly reduce the
number combinatorial possibilities for an exceptional cluster.

\subhead \refz 3.3. Contexts\endsubhead

Each upright diagonal has a context $\x(p,q)$, with $p$
the number of anchors and $p-q$ the number
of quarters around the diagonal \cite{F}.
The dihedral angle
of a quarter is less than $\pi$ ($\A_8$), so the context
$\x(2,0)$ is impossible.
There is at least
one quarter, so $p\ge q+1$, $p\ge2$.

The context $\x(2,1)$ is treated in \cite{F}.  Proposition F.4.7
shows that by removing the upright diagonal, and scoring the surrounding
region by a truncated version of the Voronoi function, an upper
bound on the score is obtained.
In the remaining contexts, $p\ge3$.
We start with contexts satisfying $p=3$.
The context $\x(3,0)$ is
to be regarded as two quasi-regular tetrahedra sharing a face rather
than as three quarters along a diagonal.  In particular, by \cite{F},
the upright quarters do not belong to the $Q$-system.

We recall that the score of an upright quarter is given by
$$\sigma(Q)=\nu(Q)=(\mu(Q)+\mu(\hat Q) + \vor_0(Q)-\vor_0(\hat Q))/2,$$
except in the contexts $\x(2,1)$ and $\x(4,0)$.  The context $\x(2,1)$
has been treated, and the context $\x(4,0)$ does
not occur in exceptional clusters.  Thus, for the remainder of this
paper, the scoring rule $\sigma(Q)=\nu(Q)$ will be used.

We have several different variants on the score depending on the
truncation, analytic continuation, and so forth.  If
$f$ is any of the functions
$$\vor_0,\ \vor,\ \Gamma, \ \nu,$$
we set $\tau_{0}$, $\tau_V$, $\tau_\Gamma$,
$\tau_\nu$, respectively,
to $$\tau_* = -f(S) +\sol(S)\zeta\pt.$$
We set $\tau(S,t) = -\vor(S,t)+\sol(S)\zeta\pt$.
The family of functions $\tau_*$ measure what is squandered by
a simplex.  We say that $Q$ has
{\it compression type\/} or {\it Voronoi type\/}
according to the scoring of $\mu(Q)$.

Crowns and anchor correction terms are used in \cite{F} to
erase upright quarters.  We imitate those methods here.
The functions $\cro$ and $\anc$ are defined and discussed
in Section F.4.
If $S=S(y_1,\ldots,y_6)$ is a simplex along $(0,v)$, set
$$\kappa(S(y_1,\ldots,y_6))=\cro(y_1/2)\dih(S)/(2\pi) +
    \anc(y_1,y_2,y_6)+\anc(y_1,y_3,y_5).$$
 $\kappa(S)$ is a bound on the difference in the
score resulting from truncation around $v$.
Assume that $S$ is the simplex formed by $(0,v)$ and two consecutive anchors
around $(0,v)$. Assume further that the circumradius of $S$ is at least
$\eta_0(y_1/2)$.  Then we have
$$\kappa(S) = -4\doct \Vol(\delta_P(W^e)),$$
where $W^e$ is the extended wedge constructed in Section \refz
2.3. To see this, it is a matter of interpreting the terms in
$\kappa$. The function {\it crown\/} enters the volume through the
region over the spherical cap $D_0$ of Section \refz 2.3, lying
outside $B(t_0)$.  By multiplying by $\dih(S)/(2\pi)$, we select
the part of the spherical cap over the unextended wedge $W$
between the anchors.  The terms {\it anc\/} adjust for the four
Rogers simplices lying above the extension $W^e$.

\subhead \refz 3.4. Three anchors\endsubhead
\smallskip

\proclaim{Lemma \refz 3.4.1}  The upright diagonal can be erased
in the context $\x(3,2)$.
\endproclaim

\demo{Proof}
Let $v_1$ and $v_2$ be the two anchors of the upright
diagonal $(0,v)$ along the quarter. Let the third
anchor be $v_3$.

Assume first that $|v|\ge 2.696$. If $Q$ is of compression type,
then by $\A_{10}$, the score is dominated by the truncated Voronoi
function $\vor_0$.  Assume $Q$ is of Voronoi type. If $|v_1|$,
$|v_2|\le 2.45$, then $\A_{11}$ gives the result. Take $|v_2|\ge
2.45$.  By symmetry, $|v-v_1|$ or $|v-v_2|\ge 2.45$. The case
$|v-v_1|\ge2.45$ is treated by $\A_{11}$. We take
$|v-v_2|\ge2.45$. Let $S=(0,v,v_2,v_3)$. If $S$ is of type $S_C$,
the result follows from $\A_{11}$. ($S$ is of type $S_C$, iff
$y_4\le 2.77$, (because $\eta_{456}\ge\eta(2.45,2,2.77)>\sqrt2$.)
If $S$ is  not of type $S_C$, we argue as follows. The function
$h^2(\eta(2h,2.45,2.45)^{-2}-\eta_0(h)^{-2})$ is a quadratic
polynomial in $h^2$ with negative values for
$2h\in[2.696,2\sqrt{2}]$. From this we find
$$\rad(S)\ge \eta(2h,2.45,2.45) \ge \eta_0(h),\quad\text{where } 2h=|v|,$$
and this justifies the use of $\kappa$ (see Section 2.3(2)).
That the truncated Voronoi function dominates the score now follows
from $\A_9$.

Now assume that $|v|\le 2.696$. If the simplices $(0,v,v_1,v_3)$
and $(0,v,v_2,v_3)$ are of type $S_C$, the bound follows from
$\A_{10},\A_{11}$. If say  $S=(0,v,v_2,v_3)$ is not of type $S_C$,
then
$$\rad(S)\ge\sqrt2>  \eta_0(2.696/2)\ge\eta_0(h),$$
justifying the use of $\kappa$. The bound follows from
$\A_9,\A_{10},\A_{11}$. ($\Gamma+\kappa <\octavor_0$,
etc.)\qed\enddemo

\smallskip

\proclaim{Lemma \refz 3.4.2}  The upright diagonal can be erased in
the context $\x(3,1)$,
 provided the three anchors
do not form a flat quarter at the origin.
\endproclaim

\demo{Proof}  In the absence of a flat quarter, truncate,
score,
and remove the vertex $v$ as in the context $\x(3,1)$ of Proposition F.4.7.  If
there is a flat quarter, by the rules of \cite{F}, $v$ is enclosed
over the flat quarter.
We do nothing further with them. This unerased case
appears in the
summary at the end of the section (3.10).\qed
\enddemo

\subhead \refz 3.5.  Six anchors\endsubhead

\proclaim{Lemma}  An upright diagonal has at most five anchors.
\endproclaim

\demo{Proof}  The proof relies on constants
and inequalities from $\A_3$ and $\A_8$.
If between two anchors there is a quarter, then
the angle is greater than $0.956$, but if there
is not,  the angle is greater than $1.23$.  So if
there are $k$ quarters and at least six anchors, they squander more than
    $$ k (1.01104) - [2\pi-(6-k)1.23]0.78701 > \squander,$$
for $k\ge0$.
\qed
\enddemo

\subhead \refz 3.6. Anchored simplices\endsubhead

Let $(0,v)$ be an upright diagonal, and let $v_1,v_2,\ldots,v_k=v_1$
be its anchors, ordered cyclically around $(0,v)$.  This cyclic
order gives dihedral angles between consecutive
anchors around the upright diagonal.  We define the dihedral angles
so that their sum is $2\pi$, even though
this will lead us to depart from our usual conventions by
assigning a dihedral angle greater than $\pi$ when all
the anchors are concentrated in some half-space bounded by
a plane through $(0,v)$.
When the dihedral angle of $S=(0,v,v_i,v_{i+1})$ is at most $\pi$,
we say that
$S$ is an {\it anchored simplex\/} if $|v_i-v_{i+1}|\le3.2$.
(The constant $3.2$ appears throughout this paper.)
All upright quarters are anchored simplices.
If an upright diagonal is completely surrounded by anchored
simplices, the configuration of anchored simplices
is sometimes called a {\it loop}.
If $|v_i-v_{i+1}|>3.2$ and the angle is less than $\pi$,
we say there is a {\it large gap\/}
around $(0,v)$ between $v_i$ and $v_{i+1}$.

To understand how anchored simplices overlap we need a bound
satisfied by vertices enclosed over an anchored simplex.

\proclaim{Lemma}  A vertex $w$ of height between 2 and $2\sqrt{2}$,
enclosed in the cone over an anchored simplex $(0,v,v_1,v_2)$ with
diagonal $(0,v)$,
satisfies $|w-v|\le 2.51$. In particular, if $|w|\le 2.51$,
then $w$ is an anchor.
\endproclaim

\demo{Proof}  As in Lemma I.3.5, the vertex $w$ cannot lie inside
the anchored simplex.  If $|v_1-v_2|\le 2\sqrt{2}$, the result
follows from Lemma F.2.2 (or Lemma F.1.9).
In fact, if $|w|\le 2\sqrt{2}$, the Voronoi
cells at $0$ and $w$ meet, so that Lemma F.2.2 forces
$(0,v_1,v_2,w)$ to be a quarter. (This observation gives a second
proof of F.1.9.)

Assume that a figure exists with $|v_1-v_2|>2\sqrt{2}$.
Suppose for a contradiction that
$|v-w|>2.51$.    Pivot $v_1$
around $(0,v_2)$ until $|v-v_1|=2.51$ and $v_2$ around $(0,v_1)$
until $|v-v_2|=2.51$.  Rescale $w$ so that $|w|=2\sqrt{2}$.
Set $x = |v_1-v_2|$.
If, through geometric considerations, $w$ is not deformed into
the plane of $(0,v_2,v_1)$, then we are left with the one-dimensional
family
$|w'|=|w'-w|=2$, for $w'=v_2,v_1$, $|v-w|=|v|=|v_1-v|=|v_2-v|=2.51$,
depending on  $x$.
This gives a contradiction
$$\align
\pi &\ge \dih(v_2,v_1,0,v) + \dih(v_2,v_1,v,w)\\
        &= 2\dih(S(x,2,2.51,2.51,2.51,2))
         > \pi,
\endalign$$
for $x>2\sqrt{2}$. (Equality is attained if $x=2\sqrt{2}$.)

Thus, we may assume that $w$ lies in the plane $P=(0,v_1,v_2)$.
Take the circle in $P$ at distance $2.51$ from $v$.
The vertices $0$ and $w$ lie on or outside the circle.
The vertices
$v_1$ and $v_2$ lie on the circle, so the diameter is at least
$x>2\sqrt{2}$.  The distance from $v$ to
$P$ is less than $x_0= \sqrt{2.51^2-2}$.  The edge $(0,w)$ cannot
pass through the center of the circle, because
$|w|$ is less than the diameter.
Reflect $v$ through $P$ to get $v'$.  Then $|v-v'|< 2x_0$.
Swapping $v_1$ and $v_2$ as necessary, we may assume that $w$ is
enclosed over $(0,v,v',v_2)$.  The desired bound $|v-w|\le 2.51$
now follows from geometric considerations and the contradiction
$$2\sqrt{2}= |w| > {\Cal E}(S(2,2.51,2.51,2x_0,2.51,2.51),2,2.51,
        2.51)= 2\sqrt{2}.$$\qed
\enddemo

\proclaim{Corollary}
A vertex of height at most $2.51$ is never enclosed
over an anchored simplex.
\endproclaim

\demo{Proof}  If so, it would be an anchor to the upright diagonal,
contrary to the assumption that the anchored simplex is formed by
consecutive anchors.\qed
\enddemo

\subhead \refz 3.7 Surrounded upright diagonals\endsubhead

This proposition is a consequence of the two lemmas that follow.
The context of the proposition is the set of anchored simplices
that have not been erased by previous reductions.

\proclaim{Proposition}  Anchored simplices do not overlap.
\endproclaim

The remaining contexts have four or  five anchors.
Let $w$ and the anchored simplex $S=(0,v,v_1,v_2)$
be as in Section \refz 3.6.
Our object is to describe the local geometry when an
upright diagonal is enclosed over an anchored simplex.
If $|v_1-v_2|\le 2\sqrt{2}$, we have seen in Section F.1.8
that there can be no enclosed upright diagonal
with $\ge 4$ anchors over the anchored
simplex $S$.

Assume  $|v_1-v_2|>2\sqrt{2}$.
Let $w_1,\ldots,w_k$, $k\ge4$,
be the anchors of $(0,w)$, indexed consecutively.
The anchors of $(0,w)$
do not lie in $C(S)$, and the
triangles $(0,w,w_i)$ and $(0,v,v_j)$ do not overlap.
Thus, the plane $(0,v_1,v_2)$ separates $w$ from $\{w_1,\ldots,w_k\}$.
Set $S_i=(0,w,w_i,w_{i+1})$. By $\A_8$,
$$\pi\ge \dih(S_1)+\cdots+\dih(S_{k-1})\ge (k-1)0.956.$$

Thus, $k=4$.
The configuration of three simplices $\{S_i\}$, which we denote
by $\Sminus$,
 will be studied
in the next two lemmas.
The superscript reminds us that $\sum\dih(S_i)-\pi$ is negative.

 We claim that $\{v_1,v_2\}=\{w_1,w_4\}$.
Suppose to the contrary that, after reindexing as necessary,
$S_0=(0,w,w_1,v_1)$ is a simplex, with $v_1\ne w_1$, that does
not overlap $S_1,\ldots,S_3$.
Then $\pi\ge \dih(S_0)+\cdots+\dih(S_3)$.
So $0.28\ge \pi-3(0.956)\ge \dih(S_0)$.
$\A_8$ now implies that $|w-v_1|\ge 2\sqrt{2}$.

Assume that $(0,w,v_1,v_2)$ are coplanar.  Disregard the
other vertices. We minimize
$|v_1-v_2|$ when
    $$|w|=2\sqrt{2},\quad |v_2|=|v_1|=|w-v_2|=2,\quad |w-v_1|=2\sqrt{2}.$$
This implies
$3.2\ge|v_1-v_2|\ge x$, where $x$ is the largest positive root
of the polynomial $\Delta(8,4,4,x^2,4,8)$.
But $x\approx 3.36$, a contradiction.

Since $(0,w,v_1,v_2)$ cannot be coplanar vertices, geometric
considerations apply and
$$2\sqrt{2}\ge |w| \ge {\Cal E}(S(2,2,2,2,2,3.2),2\sqrt{2},2,2)>2\sqrt{2}.$$
This contradiction establishes that $v_1=w_1$.

\proclaim{Lemma \refz 3.7.1}  If there is an upright diagonal $(0,v)$ with
four anchors all concentrated in a half-space through $(0,v)$, then
the three anchored simplices squander more than $0.5606$ and score
at most $-0.4339$.
\endproclaim

\demo{Proof}  The proof makes use of constants
and inequalities from $\A_2$, $\A_8$,
and $\A_{12}$.  The dihedral angles are at most
$\pi- 2(0.956) < 1.23$.
This forces $y_4\le 2.51$, for each simplex $S$.
  So they are all quarters.
The
three anchored simplices squander at least
$$3 (1.01104) - \pi (0.78701) > 0.5606.
$$
The bound on score follows similarly from $\nu<-0.9871+0.80449\dih$.
\qed\enddemo

\proclaim{Lemma \refz 3.7.2}  If an $\Sminus$ configuration
overlaps an anchored simplex,
the decomposition star squanders at least $\squander$.
\endproclaim

\demo{Proof}
Suppose that $(0,v,v_1,v_2)$ is an anchored
simplex that another anchored simplex overlaps, with $(0,v)$ the
upright diagonal.  Let $(0,w)$ be the upright diagonal of an
$\Sminus$ configuration.
We score the two simplices $S'_i = (0,v,w,v_i)$ by truncation at
$\sqrt{2}$.  Truncation at $\sqrt{2}$
is justified by face-orientation arguments or by geometric
considerations:
    $$\Cal E(S(2,2.51,2.51,2.51,2.51,2.51),2,2,2)>2\sqrt{2}.$$
By $\A_{12}$,
    $$\tau_V(S'_1,\sqrt{2})+\tau_V(S'_2,\sqrt{2})\ge 2(0.13) +
        0.2(\dih(S'_1)+\dih(S'_2)-\pi) > 0.26.$$
Together with the three simplices in $\Sminus$
 that squander
at least $0.5606$, we obtain the stated bound.\qed
\enddemo

\subhead \refz 3.8.  Five anchors\endsubhead

When there are five anchors of an upright diagonal, each dihedral
angle around the diagonal is at most $2\pi-4(0.956)<\pi$.
There are at most two large gaps by $\A_8$,
$$3(1.65)+2(0.956)>2\pi.$$

\proclaim{Lemma 3.8.1}  If an upright diagonal has five anchors with
two large gaps, then the three anchored simplices squander
$>\squander$.
\endproclaim

\demo{Proof} By $\A_8$, the anchored simplices are all quarters,
 $1.23+2(1.65)+2(0.956)>2\pi$.  The dihedral angle
is less than $2\pi-2(1.65)$.  The linear programming bound from
$\A_3$ is greater than $0.859>\squander$.
\qed
\enddemo

Define a {\it masked\/} flat quarter to be a flat quarter that is not
in the $Q$-system because it overlaps an upright quarter in
the $Q$-system.
They can only occur in a very special setting.

\proclaim{Lemma 3.8.2} Let $(0,v)$ be an upright diagonal with at least four
anchors.  If $Q$ is a flat quarter that overlaps an anchored
simplex that lies along $(0,v)$,
then the vertices of $Q$ are the origin
and three consecutive anchors of $(0,v)$.
\endproclaim

\demo{Proof}  For there to be overlap, the diagonal $(w_1,w_2)$
of $Q$ must pass through the face $(0,v,v_1)$ formed by
some anchor $v_1$.  (see Lemma F.1.3).  By Lemma
F.1.5, $w_1$ and $w_2$ are
anchors of $(0,v)$.  By Lemma F.1.8, $w_2,v_1$, and $w_1$ are consecutive anchors.
If $v_1$ is a vertex of $Q$ we are done.
Otherwise, let $w_3\ne 0,w_1,w_2$ be the remaining vertex of $Q$.
The edges $(v,v_1)$
and $(v_1,0)$ do not pass through the face $(w_1,w_2,w_3)$ by
Lemma F.1.3.
Likewise, the edges $(w_2,w_3)$ and $(w_3,w_1)$ do not pass through the
face $(0,v,v_1)$.  Thus, $v$ is enclosed over the quarter $Q$.

\was{Since $v$ has at least
four anchors it is not part of an isolated pair (Diagram F.1.10).
Thus there is another quarter $Q'$ adjacent to $Q$ along $(w_1,w_2)$.
If this quarter shares the face $(0,w_1,w_2)$ with $Q$, we have a
quad cluster, and the four anchors give us an octahedron
and the result.
If the quarter $Q'$ shares the face $(w_1,w_2,w_3)$ with $Q$,
then by Lemma F.1.8, the
other vertex of $Q'$ is $v$.
Thus, $w_3$ is an anchor of $v$.}

Let $w_3'\ne w_1,v_1,w_2$ be a fourth anchor of $(0,v)$.
By Lemma F.1.3, we have $w_3'=w_3$.
\qed
\enddemo

\proclaim{Corollary (of proof)}  If $v$ is enclosed over a
flat quarter, then $(0,v)$ has at most four anchors.\qed
\endproclaim

When we are unable to erase the upright diagonal with
five anchors and a
large gap, we are able to obtain strong bounds on
the score.  We let $\Sfour$ denote the configuration of four
upright quarters and the large gap around an upright diagonal.

\proclaim{Lemma \refz 3.8.3}  Suppose an upright diagonal has five
anchors and one large gap. The four anchored simplices score at
most $-0.25$. The four anchored simplices squander at least $0.4$.
If any of the four anchored simplices is not an upright quarter
then the four simplices squander at least $\squander$.
\endproclaim

\demo{Proof} $\A_2$ and  $\dih>1.65$ from $\A_8$ give the bound
$-0.25$. $\A_3$ gives the bound $0.4$.  To get the final statement
of the lemma, use inequalities $\A_5$ and $\A_7$ as well.
    \qed
\enddemo

\proclaim{Corollary}  There is at most one $\Sfour$.
\endproclaim

\demo{Proof}  The crown along the large gap, with the
bound of the lemma, gives
    $0.4-\kappa \ge 0.4+0.02274$
squandered by each $\Sfour$ (see $\A_9$).  The rest squanders
a positive amount (see Lemma 4.1).  If there are two $\Sfour$-configurations,
use $2(0.4+0.02274)>\squander$.\qed\enddemo


We set $\xiG = 0.01561$, $\xiV = 0.003521$, $\xiG'=0.00935$,
$\xik=-0.029$,
$\xikG = \xik+\xiG = -0.01339$.
The first two constants appear in $\A_{10}$ and $\A_{11}$ as penalties
for erasing upright quarters of compression type, and Voronoi type,
respectively.
$\xiG'$ is an improved bound on the penalty for erasing when
the upright diagonal is at least $2.57$.
Also, $\xik$ is an upper bound on $\kappa$ from $\A_9$, when
the upright diagonal is at most $2.57$.  If the upright
diagonal is at least $2.57$, then we still obtain the bound
$\xikG = -0.02274+\xiG'$ from $\A_9$ on the sum of $\kappa$ with
the penalty from erasing an upright quarter.

\subhead \refz 3.9.  Four anchors\endsubhead

\proclaim{Lemma \refz 3.9.1}
If there are at least two large gaps around an upright diagonal
with $4$ anchors, then it can be erased.
\endproclaim

\demo{Proof}
There are at least as many large gaps as upright quarters.
Each large gap drops us by $\xik$ and each quarter
lifts us by at most $\xiG$ by $\A_9,\A_{10},\A_{11}$.
We have $\xikG<0$.
\qed
\enddemo
\smallskip

\proclaim{Remark} Let $(0,v)$ be an enclosed vertex over a flat quarter.
Then $$|v|\ge\Cal E(2,2,2,2.51,2.51,2\sqrt2,2,2,2)>2.6.$$  If an
edge of the flat quarter is sufficiently short, say $y_6\le2.2$, then
$$|v|\ge\Cal E(2,2,2,2.2,2.51,2\sqrt2,2,2,2)>2.7.$$
The two dihedral angles on the gaps are $>1.65$.  If the two
quarters mask a flat quarter, we use the scoring of 3.10.2.c.
We have $0.0114< -2\xikG$.
\endproclaim

\smallskip
When there is one large gap, we may erase with a penalty $\pi_0=0.008$.

\proclaim{Lemma \refz 3.9.2} Let $v$ be an upright diagonal with 4
anchors.  Assume that there is one large gap.  The anchored
simplices can be erased with penalty $\pi_0=0.008$. If any of the
anchored simplices around $v$ is not an upright quarter then we
can erase with penalty $\pi_0=0.00222$.

Moreover, if there is a flat quarter
overlapping an upright quarter, then (1) or
(2) holds.

(1)  The truncated Voronoi function exceeds the score by at least
$0.0063$.  The diagonal of the flat is at least 2.6, and the edge
opposite the diagonal is at least 2.2.

(2)  The truncated Voronoi function exceeds the score by at least
$0.0114$.  The diagonal of the flat is at least 2.7, and the edge
opposite the diagonal is at most 2.2.
\endproclaim

As a matter of notation, we let $\Splus$ be the configuration
of three simplices described by the lemma, when there is no masked flat
quarter.

\demo{Proof}
The constants and inequalities used in this proof can be found in $\A_9$,
$\A_{10}$, and $\A_{11}$.

First we establish the penalty $0.008$.
The truncated Voronoi function is an upper bound on the score
of an anchored simplex that is not a quarter.
 By these inequalities, the result
follows  if the diagonal satisfies
$y_1\ge 2.57$.
Take $y_1\le 2.57$.

If any of the upright quarters are of Voronoi type, the result
follows from $(\xikG+\xiG<0.008)$.
If the edges along the large gap are less than 2.25, the result
follows from $(-0.03883+3\xiG = 0.008)$.
If all but one edge along the large gap are  less than 2.25, the
result follows from
$(-0.0325 + 2\xiG + 0.00928 = 0.008)$.

If there are at least two edges along the large gap of length at least
2.25, we consider two cases according to whether they lie on a
common face of an upright quarter.  The same group of inequalities
from the appendix gives the result.
The bound 0.008 is now fully established.

\smallskip
Next we prove that we can erase with penalty $0.00222$, when one
of the anchored simplices is not a quarter.  If $|v|\ge2.57$, then
we use
    $$2\xiG + \xiV + \xik \le 0.00935+0.003521 -0.2274\le 0.$$
If $|v|\le2.57$, we use
    $$2(0.01561)-0.029 \le 0.00222.$$

\smallskip
Let $v_1\ldots,v_4$ be the consecutive anchors of
the upright diagonal $(0,v)$ with $(v_1,v_4)$ the large gap.
Suppose $|v_1-v_3|\le 2\sqrt{2}$.

We claim the upright diagonal $(0,v)$ is not enclosed over
$(0,v_1,v_2,v_3)$.
Assume the contrary.  The edge $(v_1,v_3)$ passes through the
face $(0,v,v_4)$.  Disregarding the vertex $v_2$, by geometric
considerations, we arrive at the rigid figure
$$
\align
|v|&=2\sqrt{2},\ |v_1|=|v_1-v|=|v-v_3|=|v_3|=|v_3-v_4|=2
\\ |v-v_4|&=|v_4|=2.51,\ |v_1-v_4|=3.2.
\endalign
$$
The dihedral angles of $(0,v,v_1,v_4)$ and $(0,v,v_3,v_4)$ are
$$\dih(S(2\sqrt{2},2,2.51,3.2,2.51,2))>2.3,\
\dih(S(2\sqrt{2},2,2.51,2,2.51,2))>1.16\
$$
The sum is greater than $\pi$, contrary to the claim that the edge
$(v_1,v_3)$ passes through the face $(0,v,v_4)$.
(This particular conclusion leads to the corollary cited at the
end of the proof.)
Thus, $(v_1,v_3)$ passes through $(0,v,v_2)$ so that the
simplices $(0,v,v_1,v_2)$
and $(0,v,v_2,v_3)$ are of Voronoi type.

To complete the proof of the lemma, we show that when
there is a masked flat quarter, either (1) or (2) holds.
Suppose we mask a flat quarter $Q'=(0,v_1,v_2,v_3)$.
We have established that $(v_1,v_3)$ passes through
the face $(0,v,v_2)$.
To establish (1) assume that $|v_2|\ge 2.2$.  The remark
before the lemma gives
$$|v_1-v_3|\ge \Cal E(S(2,2,2,2\sqrt{2},2.51,2.51),2,2,2)>2.6.$$
The bound $0.0063$ comes from
$$\xikG + 2\xiV < -0.0063$$

To establish (2) assume that $|v_2|\le 2.2$.
The remark gives
$$|v_1-v_3|\ge \Cal E(S(2,2,2,2\sqrt{2},2.2,2.51),2,2,2)>2.7.$$
  If the simplex
$(0,v,v_3,v_4)$ is of Voronoi type, then
$$\xik + 3\xiV  < -0.0114$$
Assume that $(0,v,v_3,v_4)$ is of compression type.
We have
$$-0.004131 +\xikG + \xiV \le -0.0114.$$\qed
\enddemo

\proclaim{Corollary (of proof)}  If there are four anchors and if
the upright diagonal is enclosed over a flat quarter, then there
are four anchored simplices and at least three quarters around the
upright diagonal.
\qed
\endproclaim

\subhead \refz 3.10. Summary\endsubhead

The following index summarizes the cases of upright quarters
that have been treated in Section \refz 3.
If the number of anchors is the number of anchored simplices (no
large gaps), the results appear in Section 5.11.
Every other possibility has been treated.

\bul {0,1,2 anchors\hfill Sec. \refz 3.3}

\bul 3 anchors \hfill Sec. \refz 3.4

\quad\cir context $\x(3,0)$

\quad\cir context $\x(3,1)$

\quad\cir context $\x(3,2)$

\quad\cir context $\x(3,3)$

\bul 4 anchors \hfill Sec. \refz 3.9

\quad \cir 0 gaps (Section 5.11)

\quad\cir 1 gap

\quad\cir 2 or more gaps

\bul 5 anchors \hfill Sec. \refz 3.8

\quad\cir 0 gaps (Section 5.11)

\quad\cir  1 gap ($\Sfour$)

\quad\cir 2 or more gaps

\bul {6 or more anchors \hfill Sec. \refz 3.5}

\smallskip
By truncation and various comparison lemmas,
we have entirely eliminated upright diagonals except when
there are between three and five anchors.  We may assume that there is at most one
large gap around the upright diagonal.

\smallskip
\reft 1.  Consider an anchored simplex $Q$ around a remaining upright diagonal.
The score of is $\nu(Q)$ if $Q$ is a quarter, the analytic
Voronoi function $\vor(Q)$ if the
simplex is of type $S_C$ (Section \refz 2.5),
and the truncated Voronoi function $\vor_0(Q)$ otherwise.

\smallskip
\reft 2.  Consider a flat quarter $Q$
in an exceptional cluster.  An upper bound
on the score is obtained by taking the maximum
of all of the following
functions that satisfy the stated conditions on $Q$.  Let $y_4$ denote
the length of the diagonal and $y_1$ be the length of the opposite edge.

(a)  The function $\mu(Q)$.

(b)  $\vor_0(Q) - 0.0063$, if $y_4\ge 2.6$ and $y_1\ge 2.2$.\hfill
    (Lemma \refz 3.9)

(c)  $\vor_0(Q) - 0.0114$, if $y_4\ge 2.7$ and $y_1\le 2.2$.
    \hfill (Lemma \refz 3.9)

(d)  $\nu(Q_1)+\nu(Q_2)+\vor_x(S)$, if there is an enclosed vertex
    $v$ over $Q$ of height between $2.51$ and $2\sqrt{2}$ that
    partitions the convex hull of $(Q,v)$ into two upright quarters
    $Q_1$, $Q_2$ and a third simplex $S$. Here $\vor_x=\vor$
    if $S$ is of type $S_C$, and $\vor_x=\vor_0$ otherwise.
    \hfill (Lemma \refz 3.4)

(e)  $\vor(Q,1.385)$ if the simplex is of type $S_B$ (Section \refz 2.5).

(f) $\vor_0(Q)$ if the simplex is an isolated quarter with
    $\max(y_2,y_3)\ge2.23$, $y_4\ge2.77$,
    and $\eta_{456}\ge\sqrt2$.

\smallskip

\reft 3.   If $S$ is a simplex is of type $S_A$, its score is $\vor(S)$.
(Section \refz 2.5.)

\smallskip

\reft 4.  Everything else is scored by the truncation of Voronoi, $\vor_0$.
    The Formula F.3.7 is used on these remaining pieces.
    On top of what is obtained for the standard cluster by summing all
these terms,
there is a penalty $\pi_0=0.008$ each time the simplex configuration
$\Splus$ is erased.

\smallskip
\reft 5.  The remaining upright diagonals not surrounded by
    anchored simplices are
    $\Splus$, $\Sminus$,  $\Sfour$
    from Section 3.7, 3.8 and 3.9.

\bigskip
\subhead{\refz 3.11. Some flat quarters}\endsubhead 

Recall that $\xiV=0.003521$, $\xiG=0.01561$, $\xiG'=0.00935$. They
are the penalties that result from erasing an upright quarter of
Voronoi type, an upright quarter of compression type, and an
upright quarter of compression type with diagonal $\ge2.57$.
  (See $\A_{10}$ and  $\A_{11}$.)

In the next lemma, we score a flat quarter by any of the functions
on the given domains
     $$\hat\sigma= \cases \Gamma,& \eta_{234},\eta_{456}\le\sqrt2,\\
             \vor, &\eta_{234}\ge\sqrt2,\\
            \vor_0, & y_4\ge 2.6, y_1\ge2.2,\\
            \vor_0, & y_4\ge 2.7,\\
            \vor_0,& \eta_{456}\ge\sqrt2.
    \endcases
    $$

\proclaim{Lemma \refz 3.11.1}  $\hat\sigma$ is an upper bound on
the functions in Section 3.10(a)--(f). That is, each function in
Section 3.10 is dominated by some choice of $\hat\sigma$.
\endproclaim

\demo{Proof}  The only case in doubt is the function of 3.10(d):
$$\nu(Q_1)+\nu(Q_2)+\vor_x(S).$$ This is established by the
following lemma.
\qed
\enddemo

We consider the context $\x(3,1)$ that occurs when two upright
quarters in the $Q$-system lie over a flat quarter. Let $(0,v)$ be
the upright diagonal, and assume that $(0,v_1,v_2,v_3)$ is the
flat quarter, with diagonal $(v_2,v_3)$. Let $\sigma$ denote the
score of the upright quarters and other anchored simplex lying
over the flat quarter.

\proclaim{Lemma \refz 3.11.2} $\sigma\le \min(0,\vor_0)$.
\endproclaim

\demo{Proof}  The bound of $0$ is established in \cite{II} and
\cite{F}.

By F.4.7.5, if $|v|\ge 2.69$, then the upright quarters satisfy
$$\nu < \vor_0 + 0.01 (\pi/2-\dih)$$
so the upright quarters can be erased.  Thus we assume without
loss of generality that $|v|\le 2.69$.

We have $$|v|\ge\E(S(2,2,2,2t_0 ,2t_0 ,2\sqrt2),2,2,2)>2.6.$$ If
$|v_1-v_2|\le 2.1$,  or $|v_1-v_3|\le 2.1$, then
    $$|v|\ge \E(S(2,2,2,2.1,2t_0 ,2\sqrt2),2,2,2)>2.72,$$
contrary to assumption.  So take $|v_1-v_2|\ge 2.1$ and
$|v_1-v_3|\ge2.1$. Under these conditions we have the interval
calculation $\nu(Q) < \vor_0(Q)$ where $Q$ is the upright quarter
(see $\A_{13}$). \qed\enddemo

\proclaim{Remark}  If we have an upright diagonal enclosed over a
masked flat quarter in the context $\x(4,1)$, then there are 3
upright quarters.  By the same argument as in the lemma, the two
quarters over the masked flat quarter score $\le\vor_0$. The third
quarter can be erased with penalty $\xiV$.
\endproclaim

Define the {\it central vertex\/} $v$ of a flat quarter to be the
vertex for which $(0,v)$ is the edge opposite the diagonal.

\proclaim{Lemma \refz 3.11.3} $\mu < \vor_0 +0.0268$ for all flat
quarters. If the central vertex has height $\le2.17$, then
$\mu<\vor_0+0.02$.
\endproclaim

\demo{Proof} This is an interval calculation.  See
$\A_{13}$.\qed\enddemo

We  measure what is squandered by a flat quarter by $\hat\tau =
\sol\zeta\pt - \hat\sigma$.

\proclaim{Lemma \refz 3.11.4}  Let $v$ be a corner of an
exceptional cluster at which the dihedral angle is at most $1.32$.
Then the vertex $v$ is the central vertex of a flat quarter $Q$ in
the exceptional region. Moreover, $\hat\tau(Q)>3.07\,\pt$. If
$\hat\sigma=\vor_0$ (and if $\eta_{456}\ge\sqrt2$), we may use the
stronger constant $\tau_0(Q)> 3.07\,\pt+\xi_V+2\xiG'$.
\endproclaim

\demo{Proof} Let $S=S(y_1,\ldots,y_6)$ be the simplex inside the
exceptional cluster centered at $v$, with $y_1=|v|$. The
inequality $\dih\le 1.32$ gives the interval calculation $y_4\le
2\sqrt{2}$, so $S$ is a quarter. The result now follows by
interval arithmetic. See $\A_{13}$.\qed\enddemo


\heads{\refz 4. Distinguished Edges and Subregions}
\vfill\eject
\head \refz 4.  Distinguished Edges and Subregions\endhead

\subhead \refz 4.1. Positivity\endsubhead

\proclaim{Lemma}  $\tau_{0}\ge 0$ on local $V$-cells.
\endproclaim

\demo{Proof}  Everything truncated at $t_0$ can be broken into
three types of pieces: Rogers simplices $R(a,b,t_0)$, wedges of
$t_0$-cones,
and spherical regions. (See Diagram II.4.2.)
The wedges of $t_0$-cones and spherical regions can be considered as the
degenerate cases $b=t_0$ and $a=b=t_0$ of Rogers simplices, so it
is enough to show that $\tau(R(a,b,t_0))\ge 0$.
We have $t_0>\sqrt{3/2}$, so by Rogers's lemma (I.8.6.2),
$$\tau(R(a,b,t_0))>\tau(R(1,\eta(2,2,2),\sqrt{3/2})).$$
The right-hand side is zero. (In fact, the vanishing
of the right-hand side is essentially Rogers's bound.  Nothing
is squandered when Rogers's bound is met.)
\qed\enddemo

\subhead \refz 4.2.  Distinguished edge conditions\endsubhead

Take an exceptional cluster.  We prepare the cluster by erasing
all upright diagonals possible, including the $\Splus$, $\Sfour$, $\Sminus$ configurations.
The only remaining upright diagonals are
those on an upright diagonal surrounded by anchored
simplices (loops).
When the upright diagonal is erased, we score
with the truncated Voronoi function.
The exceptional clusters in Sections 4 and 5 are assumed
to be prepared in this way.

A simplex $S$ is {\it special\/} if the
fourth edge has length at least $2\sqrt{2}$
and at most $3.2$,
and the others have length at most $2.51$.
The fourth edge will be called its diagonal.

We draw a system of edges between
vertices.  Each vertex will have height at most $2.51$.
  The projections of the edges
to the unit sphere will divide the standard region into subregions.
We call an edge {\it nonexternal\/} if the projection of the edge
    lies entirely in the (closed) exceptional region.

\reft 1.  Draw all nonexternal edges of length at most $2\sqrt{2}$ except
    those between nonconsecutive anchors of a remaining upright
    diagonal.
These edges do not cross (Lemma F.1.6).
These edges do not cross the edges of anchored simplices
    (Lemmas \refz 3.6, F.1.5).

\reft 2.  Draw all edges of (remaining) anchored upright simplices
that are opposite the upright diagonal,
    except when the edge gives a special simplex.
The anchored simplices do not overlap (Lemma \refz 3.7), so these
    edges do not cross.
    These edges are nonexternal (Lemma \refz 3.6, F.1.3).

\reft 3.  Draw as many additional nonexternal edges as possible of
    length at most 3.2 subject to not crossing another edge, not
crossing any edge of an anchored simplex, and not being the diagonal of a
special simplex.

We fix once and for all a maximal collection of edges subject
to these constraints.
Edges in this collection are called {\it distinguished\/} edges.
The projection of the distinguished edges to the unit sphere
gives the bounding edges of regions called the {\it subregions}.  Each
standard region is a union of subregions.
The vertices of height at most $2.51$ and the vertices of the
remaining upright diagonals are said
to form a {\it subcluster}.

By construction, the special simplices and anchored simplices
around an upright quarter form a
subcluster.  Flat quarters in the $Q$-system, flat quarters
of an isolated pair, and simplices of type $S_A$ and $S_B$ are
subclusters.
      Other
subclusters are scored by the truncation of the Voronoi function.
For these subclusters,
the Formula F.3.7 extends without modification.

\subhead \refz 4.3.  Scoring subclusters\endsubhead

The terms of Equation F.3.7 defining $\vor_0(P)=\vor(P,t_0)$ have a clear
geometric interpretation as quoins, wedges of $t_0$-cones, and solid angles (see
\cite{F}). There is a quoin for each Rogers simplex.
  There is a somewhat delicate point that arises
in connection with the geometry of subclusters.  It is not true
in general
that the Rogers simplices entering into the truncation $\vor_0(P)$
of $P$ lie in the cone over $P$.
  Formula F.3.7 should be
viewed as an analytic continuation that has a nice geometric interpretation
when things are nice, and which always gives the right answer
when summed over all the subclusters in the cluster, but which
may exhibit unusual behavior in general.
The following lemma shows that the simple geometric interpretation
of Formula F.3.7 is valid when the subregion is not triangular.

\proclaim{Lemma}
If a subregion is not a triangle and is not  the subregion containing
the anchored simplices around an upright diagonal,
 the cone of arcradius
$$\psi =\arccos(|v|/2.51)$$
centered along $(0,v)$, where $v$ is a corner of the subcluster,
does not cross out of the subregion.
\endproclaim

\demo{Proof}
For a contradiction, let $(v_1,v_2)$ be a distinguished
 edge that the cone crosses.
If both edges $(v,v_1)$ and $(v,v_2)$ have length less than
$2.51$, there can be no enclosed vertex $w$ of height at most 2.51,
unless its distance from $v_1$ and $v_2$ is less than $2.51$:
    $$\Cal E(S(2,2,2,2.51,2.51,3.2),2.51,2,2)>2.51.$$
In this case, we can replace $(v_1,v_2)$ by an edge of the subregion
closer to $v$, so without loss of generality we may assume that there
are no enclosed vertices when both edges $(v,v_1)$ and $(v,v_2)$ have
length less than $2.51$.

The subregion is not a triangle, so
$|v-v_1|\ge 2.51$, or $|v-v_2|\ge 2.51$, say $|v-v_1|\ge 2.51$.
Also $|v-v_2|\ge2$.  Pivot so that $|v_1-v_2|=3.2$, $|v-v_1|=2.51$,
$|v-v_2|=2$.  (The simplex $(0,v_1,v_2,v)$ cannot collapse ($\Delta\ne0$)
as we pivot. See Inequality 4.8(*) below.)
Then use $\beta_\psi\le \dih_3$ from $\A_1$.
\qed
\enddemo

As a consequence, in nonspecial standard regions,
the terms in the Formula F.3.7 for $\vor_0$
retain their interpretations as quoins, Rogers simplices,
$t_0$-cones, and solid angles, all lying in the cone over the standard region.
\smallskip

\subhead \refz 4.4. The main theorem\endsubhead

Let $R$ be a standard cluster. Let $U$ be the set of corners,
that is, the set of vertices in the cone over $R$ that have
height at most $2.51$.  Consider the set $E$ of
edges of length at most $2.51$ between vertices of $U$.
We attach a multiplicity to each edge.
We let the multiplicity be $2$ when the edge projects to the interior
of the standard region, and $0$ when the edge
projects to the complement of the standard region.  The other
edges, those bounding the standard region, are counted with
multiplicity $1$.

Let $n_1$ be the number of edges
in $E$,
counted with multiplicities.
Let $c$ be the number of classes of vertices
under the equivalence relation $v\sim v'$ if there is a sequence
of edges in $E$ from $v$ to $v'$.
  Let $n(R)=n_1+2(c-1)$.
If the standard region under $R$ is a polygon, then $n(R)$ is the
number of sides.

\proclaim{Theorem} $\tau(R) > t_n$, where $n=n(R)$ and
$$\align
t_4&=0.1317,\quad t_5=0.27113,\quad
t_6=0.41056,\\
t_7&=0.54999,\quad
t_8=0.6045.
\endalign$$
The decomposition star scores less than $8\,\pt$, if $n(R)\ge 9$,
for some standard cluster $R$.
The scores satisfy $\sigma(R)<s_n$, for $5\le n\le 8$, where
$$
s_5=-0.05704,\quad s_6=-0.11408,\quad
s_7=-0.17112,\quad s_8=-0.22816.
$$
\endproclaim

Sometimes, it is convenient to calculate these bounds as a multiple
of $\pt$.  We have
$$
\align
t_4&>2.378\,\pt,\quad t_5>4.896\,\pt,\quad
t_6>7.414\,\pt,\\
t_7&>9.932\,\pt,\quad
t_8>10.916\,\pt.
\endalign$$
$$
s_5 < -1.03\,\pt,\quad s_6<-2.06\,\pt,\quad
s_7<-3.09\,\pt,\quad s_8<-4.12\,\pt.
$$

\proclaim{Corollary}  Every standard region is a either a polygon or one
shown in the diagram.
\endproclaim

\smallskip
\gram|2.2||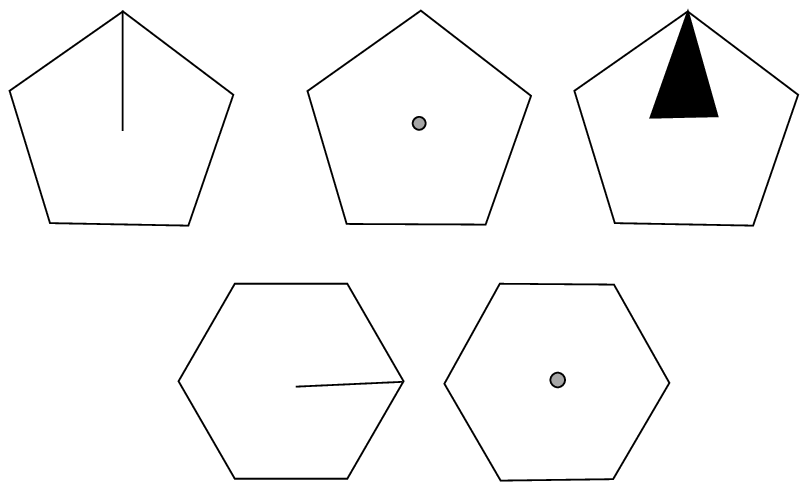|
\smallskip

In the cases that are not (simple) polygons, we call the
{\it polygonal hull\/} the
polygon obtained by removing the internal edges and vertices.
We have $m(R)\le n(R)$, where the
constant $m(R)$ is the number of sides of the polygonal hull.

\demo{Proof}  By the theorem,
if the standard region is not a polygon, then
$8\ge n_1\ge m\ge 5$. (Quad clusters and quasi-regular tetrahedra
have no enclosed vertices. See Lemma I.3.7 and Lemma III.2.2.)
If $c>1$, then $8\ge n=n_1+2(c-1)\ge 5+2(c-1)$, so $c=2$, and
$n_1=5,6$ (frames 2 and 5 of the diagram).

Now take $c=1$.    Then $8\ge n\ge 5+(n-m)$,
so $n-m\le 3$.  If $n-m=3$, we get frame 3.
If $n-m=2$, we have $8\ge m+2\ge 5+2$, so $m=5,6$ (frames 1, 4).

But $n-m=1$ cannot occur, because a single edge that does not bound
the polygonal hull has
even multiplicity.  Finally, if $n-m=0$, we have a  polygon.\qed
\enddemo

\subhead \refz 4.5. Proof\endsubhead

The proof of the theorem occupies the rest of the paper.  We begin with
a slightly simplified account of the method of proof.
Set $t_9=0.6978$, $t_{10}= 0.7891$, $t_n=\squander$,
for $n\ge 11$.
Set $D(n,k) = t_{n+k} - 0.06585\,k$, for $0\le k\le n$, and $n+k\ge 4$.
This function satisfies
$$D(n_1,k_1)+D(n_2,k_2)\ge D(n_1+n_2-2,k_1+k_2-2).\tag \refz 4.5.1$$
In fact, this inequality unwinds to $t_r+0.13943\ge t_{r+1}$,
$D(3,2)=0.13943$, and $t_n =(0.06585)2+(n-4)D(3,2)$, for $n=4,5,6,7$.
These hold  by inspection.

Call an edge between two vertices of height at most $2.51$ {\it long\/}
if it has length
greater than 2.51.
Add the distinguished edges to break the standard regions
into subregions.
We say that a subregion has {\it edge parameters}
$(n,k)$ if there
are $n$ bounding edges, where $k$ of them are long.
(We count edges with multiplicities as in Section 4.4,
if the subregion is
not a polygon.)  Combining
two subregions of edge parameters $(n_1,k_1)$ and $(n_2,k_2)$ along a long edge
$e$
gives a union with edge parameters $(n_1+n_2-2,k_1+k_2-2)$,
where we agree not to count the internal edge $e$ that no longer bounds.
  Inequality \refz 4.5.1
localizes the main theorem to what is squandered by subclusters.
Suppose we break the standard cluster into groups of subregions such
that if the group has edge parameters $(n,k)$, it squanders at least $D(n,k)$.
Then by superadditivity (\refz 4.5.1), the full standard cluster $R$ must
squander $D(n,0) = t_n$, $n=n(R)$, giving the result.

Similarly, define constants $s_4=0$,
$s_9 = -0.1972$, $s_{n}=0$,
for $n\ge10$.  Set $Z(n,k) = s_{n+k}-k\epsilon$, for $(n,k)\ne (3,1)$,
and $Z(3,1)=\epsilon$,
where $\epsilon=0.00005$ (cf. $\A_1$).
The function $Z(n,k)$ is subadditive:
    $$Z(n_1,k_1)+Z(n_2,k_2) \le Z(n_1+n_2-2,k_1+k_2-2).$$
In fact, this easily follows from $s_a+s_b\le s_{a+b-4}$, for
$a,b\ge 4$, and
$\epsilon>0$.
It will be enough in the proof of Theorem 4.4
to show that the score of a union of subregions with
edge parameters $(n,k)$ is at most $Z(n,k)$.

\subhead \refz 4.6. Nonagons\endsubhead

A few additional
comments are needed to eliminate $n=9,10$, even after the bounds
$t_9$, $t_{10}$ are established.
Suppose that $n=9$, and that $R$
squanders at least $t_9$ and
scores less than $s_9$.  This bound is already sufficient
to conclude that there are no other standard clusters except
quasi-regular tetrahedra ($t_9+t_4>\squander$).
There are no vertices of type $(4,0)$ or $(6,0)$:
$t_9+4.14\,\pt>\squander$ \cite{III.5.2}.   So all vertices not over
the exceptional cluster are of type $(5,0)$.
Suppose that there are $\ell$ vertices of type $(5,0)$.
The polygonal hull of $R$ has $m\le 9$ edges.
There are $m-2+2\ell$ quasi-regular tetrahedra.
If $\ell\le 3$,
then by III.5.3, the score is less than
    $$s_9+ (m-2+2\ell)\,\pt -0.48 \ell\,\pt < \score.$$
If on the other hand, $\ell\ge 4$, the decomposition star
squanders more than
    $$t_9+ 4(0.55)\,\pt > \squander.$$
The bound $s_9$ will be established as part of the proof of
Theorem 4.4.

The case $n=10$ is similar.  If $\ell=0$, the score is less than
    $(m-2)\,\pt\le \score$, because the score of an exceptional cluster
is strictly negative, \cite{F.3.13}.  If $\ell>0$, we squander at
least $t_{10}+ 0.55\,\pt > \squander$ (III.5.3).

\subhead \refz 4.7 Preparation of the standard cluster\endsubhead

Fix a standard cluster.  We return to the construction of subregions
and distinguished edges, to describe
 the penalties.
 Take the penalty of $0.008$ for $\Splus$.
Take the penalty $0.03344 = 3\xiG+\xikG$ for
$\Sfour$.
Take the penalty $0.04683=3\xiG$ for $\Sminus$.
Set $\maxpi=0.06688$.  The penalty in the next lemma refers
to the combined penalty from erasing all $\Sminus$, $\Splus$,
and $\Sfour$ configurations in the decomposition star.  The
upright quarters that completely surround an upright diagonal
(loops) are not erased.

\proclaim{Lemma}
The total penalty from a decomposition star is at most
$\maxpi$.
\endproclaim

 \demo{Proof}
  Before any upright quarters are erased, each
quarter squanders $>0.033$ ($\A_{13}$),
so the star squanders $>\squander$
if there are $\ge25$ quarters.  Assume there are at most $24$ quarters.
If the only penalties are $0.008$, we have $8(0.008)<\maxpi$.
If we have the penalty $0.04683$, there are at most $7$ other
quarters ($0.5606+8(0.033)>\squander$) (Lemma 3.7),
and no other penalties from this
type or from $\Sfour$, so the total penalty is at most
$2(0.008)+ 0.04683 < \maxpi$.
  Finally, if there
is one context $\Sfour$, there are at most $12$ other quarters (Section 3.8),
and erasing gives the penalty $0.03344+4 (0.008)<\maxpi$.
\qed
\enddemo

The remaining upright diagonals are surrounded by anchored simplices.
If the edge opposite the diagonal in an anchored simplex has length
$\ge2\sqrt2$, then there may be an adjacent special simplex whose
diagonal is
that edge.  Section 5.11 will give bounds on the aggregate of these
anchored simplices and special simplices.  In all other contexts,
the upright quarters have been erased with penalties.

Break the standard cluster into subclusters as in Section 4.2.
If the subregion is a triangle, we refer to the bounds of 5.7.
Sections 4.8--5.10 give bounds for subregions that are not
triangles in which all the upright quarters have been erased.
We follow the strategy outlined in Section 4.5, although the penalties
will add certain complications.

We now assume that we have a subcluster without quarters and whose
region is not triangular.  The truncated Voronoi function $\vor_0$
is the score.  Penalties are largely disregarded until Section
\refz 5.4.

We describe a series of deformations of the subcluster that increase
$\vor_0$ and decrease $\tau_0$.  These deformations disregard the
broader geometric context of the subcluster.  Consequently, we
cannot claim that the deformed subcluster exists in any decomposition
star.  As the deformation progresses, an edge $(v_1,v_2)$, not
previously distinguished, can emerge with the properties of
a distinguished edge.  If so, we add it to the collection of distinguished
edges, use it if possible to divide the subcluster into smaller
subclusters, and continue to deform the smaller pieces.  When
triangular regions are obtained, they are set aside until Section 5.7.

\subhead \refz 4.8 Reduction to polygons\endsubhead

By deformation, we can produce subregions whose boundary is
a polygon.  Let $U$ be the set of vertices over the subregion
of height $\le2.51$.  As in Section 4.4, the distinguished edges
partition $U$ into equivalence classes.  Move the vertices in
one equivalence class $U_1$ as a rigid body preserving heights until
the class comes sufficiently close to form a distinguished edge
with another subset.  Continue until all the vertices are
interconnected by paths of distinguished edges.  $\vor_0$ and
$\tau_0$ are unchanged by these deformations.

If some vertex $v$ is connected to three or more vertices by
distinguished edges, it follows from the connectedness of the open
subregion that there is more than one connected component $U_i$
(by paths of distinguished edges) of $U\setminus\{v\}$.
Move $U_1\cup \{v\}$ rigidly preserving heights and keeping
$v$ fixed until a distinguished edge forms with another component.
Continue until the distinguished edges break the subregions into
subregions with polygon boundaries.  Again $\vor_0$ and
$\tau_0$ are unchanged.

By the end of Section 4, we will deform all subregions into
convex polygons.

\proclaim{Remark}
 We will deform in such a way that the edges
$(v_1,v_2)$
will maintain a length of at least 2.
The proof that distances of at least $2$ are
maintained is given as \cite{HM,Lemma \refz 7.6}.  That proof
uses a parameter slightly larger than $2.51$, and hence it gives
a result that is slightly stronger than what is needed here.

We will deform in such a way that no vertex  crosses a boundary
of the subregion passing from outside to inside.
\endproclaim

Edge length constraints prevent a vertex from crossing a boundary
of the subregion from the inside to outside.  In fact, if $v$ is
to cross the edge $(v_1,v_2)$, the simplex $S=(0,v_1,v,v_2)$
attains volume 0.  We may assume, by the argument of the
proof of Lemma 4.3, that there are no vertices enclosed over $S$.
 Because we are assuming that the subregion
is not a triangle, we may assume that $|v-v_1|>2.51$.
We have $|v|\in[2,2.51]$.  If $v$ is to cross $(v_1,v_2)$, we may
assume that the dihedral angles of $S$ along $(0,v_1)$, and $(0,v_2)$
are acute.  Under these constraints, by the explicit formulas of
I.8, the vertex $v$ cannot cross out of the subregion
$$\Delta(S)\ge \Delta(2.51^2,4,4,3.2^2,4,2.51^2)>0.\tag*$$

We say that a corner $v_1$ is {\it visible}
from another $v_2$ if $(v_1,v_2)$
lies over the subregion.  A deformation may make $v_1$ visible
from $v_2$,
 making it a candidate for a new distinguished edge.
If $|v_1-v_2|\le 3.2$, then as soon as the deformation brings
them into visibility (obstructed until then by some $v$),
then $(*)$ shows that $|v_1-v|,|v_2-v|\le2.51$.
So $v_1,v,v_2$ are consecutive edges on the polygonal boundary,
and $|v_1-v_2|\ge 2\sqrt{4-t_0^2} > \sqrt{8}$.
By the distinguished edge conditions for special simplices,
$(v_1,v_2)$ is too long to be
distinguished.  In other words, there can be
no potentially distinguished edges hidden behind corners. They
are always formed in full view.

\subhead \refz 4.9 Some deformations\endsubhead

Consider three consecutive corners $v_3,v_1,v_2$ of a subcluster
$R$ such that the dihedral angle of $R$ at $v_1$ is greater than
$\pi$.  We call such an corner {\it concave}.  (If the angle is less
than $\pi$, we call it {\it convex}.)

Let $S=S(y_1,\ldots,y_6)=(0,v_1,v_2,v_3)$, $y_i=|v_i|$.
Suppose that $y_6>y_5$.  Let $x_i=y_i^2$.

\proclaim{Lemma \refz 4.9.1}  At a concave vertex,
$\partial \vor_0/\partial x_5 >0$ and
$\partial \tau_0/\partial x_5<0$.
\endproclaim

\demo{Proof}  As $x_5$ varies, $\dih_i(S)+\dih_i(R)$ is constant
for $i=1,2,3$.  The part of Formula F.3.7 for $\vor_0$
that depends on $x_5$ can
be written
$$-B(y_1)\dih(S)-B(y_2)\dih_2(S)-B(y_3)\dih_3(S)-4\doct
    (\quo(R_{135})+\quo(R_{315})),$$
where $B(y_i)=A(y_i/2)+\phi_0$,
$R_{135}=R(y_1/2,b,t_0)$, $R_{315}=R(y_3/2,b,t_0)$, $b=\eta(y_1,y_3,y_5)$,
and $A(h) = (1-h/t_0)(\phi(h,t_0)-\phi_0)$.
Set $u_{135}=u(x_1,x_3,x_5)$, and $\Delta_i = \partial \Delta/\partial x_i$.
(The notation comes from I.8 and F.3.)
We have
$${\partial \quo(R(a,b,c))\over \partial b} =
    {-a (c^2-b^2)^{3/2}\over 3 b (b^2-a^2)^{1/2}}\le 0$$
and $\partial b/\partial x_5\ge0$.  Also, $u\ge0$,
$\Delta\ge0$ (see I.8).  So it is enough to show
$$V_0(S) = u_{135}\Delta^{1/2}
    {\partial\over \partial x_5} (B(y_1)\dih(S)+B(y_2)\dih_2(S)
        + B(y_3)\dih_3(S))< 0.$$
By the explicit formulas of I.8, we have
$$V_0(S) = -B(y_1)y_1\Delta_6 + B(y_2)y_2 u_{135} - B(y_3)y_3 \Delta_4.$$
For $\tau_0$, we replace $B$ with $B-\zeta\pt$. It is enough to
show that
$$V_1(S) = -(B(y_1)-\zeta\pt)y_1\Delta_6 + (B(y_2)-\zeta\pt)y_2 u_{135} -
        (B(y_3)-\zeta\pt)y_3 \Delta_4<0.$$
The lemma now follows from $\A_{14}$. We note that the polynomials
$V_i$ are linear in $x_4$, and $x_6$, and this may be used to reduce
the dimension of the calculation.
\enddemo
\qed

We give a second form of the lemma when the dihedral angle of $R$ is
less than $\pi$, that is, at a convex corner.

\proclaim{Lemma \refz 4.9.2}  At a convex corner,
$\partial \vor_0/\partial x_5 <0$ and
$\partial \tau_0/\partial x_5>0$, if $y_1,y_2,y_3\in[2,2.51]$,
$\Delta\ge0$, and
(i) $y_4\in[2\sqrt{2},3.2]$, $y_5,y_6\in[2,2.51]$, or
(ii) $y_4\ge 3.2, y_5,y_6\in[2,3.2]$.
\endproclaim

\demo{Proof} We adapt the proof of the previous lemma.  Now
$\dih_i(S)-\dih_i(R)$ is constant, for $i=1,2,3$, so the signs
change.  $\vor_0$ depends on $x_5$ through
$$\sum B(y_i)\dih_i(S) - 4\doct (\quo(R_{135})+\quo(R_{315})).$$
So it is enough to show that
$$V_0 - 4\doct\Delta^{1/2}u_{135}{\partial\over\partial x_5}
    (\quo(R_{135})+\quo(R_{315}))<0.$$
Similarly, for $\tau_0$, it is enough to show that
$$V_1 - 4\doct\Delta^{1/2}u_{135}{\partial\over\partial x_5}
    (\quo(R_{135})+\quo(R_{315}))<0.$$
By $\A_{14}$
$$\align
    -4\doct  u_{135}{\partial\phantom{x}\over\partial x_5}
    (\quo(R_{135})+\quo(R_{315}))&< 0.82,\quad\hbox{on } [2,2.51]^3,\\
                            &<0.5,\quad\hbox{on }[2,2.51]^3, y_5\ge2.189.
\endalign
$$
The result now follows from
the inequalities $\A_{14}$.
\qed
\enddemo

Return to the situation of concave corner $v_1$.
Let $v_2$, $v_3$ be the adjacent corners.
By increasing $x_5$, the vertex $v_1$ moves away from every corner
$w$ for which $(v_1,w)$ lies outside the region.  This deformation
then satisfies the constraint of Remark 4.8.
Stretch the shorter
of $(v_1,v_2)$, $(v_1,v_3)$ until $|v_1-v_2|=|v_1-v_3|=3.07$
(or until a new distinguished edge forms, etc.).  Do this at
all concave corners.

By stopping at $3.07$, we prevent a corner crossing an edge from
outside-in. Let $w$ be a corner that threatens to cross a
distinguished edge $(v_1,v_2)$ as a result of the motion at a
nonconvex vertex.  To say that the crossing of the edge is from
the outside-in implies more precisely that the vertex being moved
is an endpoint, say $v_1$, of the distinguished edge.  At the
moment of crossing the simplex $(0,v_1,v_2,w)$ degenerates to a
planar arrangement, with the projection of $w$ lying over the
geodesic arc connecting the projections of $v_1$ and $v_2$.  To
see that the crossing cannot occur, it is enough to note that the
volume of a simplex with opposite edges of lengths at most $2t_0$
and $3.07$ and other edges at least $2$ cannot be planar.  The
extreme case is
    $$\Delta(2^2,2^2,(2t_0)^2,2^2,2^2,3.07^2) > 0.$$

If $|v_1|\ge2.2$, we can continue the deformations even further.
We stretch the shorter of $(v_1,v_2)$ and $(v_1,v_3)$ until
$|v_1-v_2|=|v_1-v_3|=3.2$ (or until a new distinguished edge forms,
etc.).  Do this at all concave corners $v_1$
for which $|v_1|\ge2.2$.  To see that corners cannot cross an edge
from the outside-in, we argue as in the previous paragraph,
but replacing  $3.07$
with $3.2$.  The extreme case becomes
    $$\Delta(2.2^2,2^2,(2t_0)^2,2^2,2^2,3.2^2) > 0.$$

\subhead \refz 4.10 Truncated corner cells\endsubhead

Because of the arguments in the Section \refz 4.9, we may assume
without loss of generality that we are working with a subregion
with the following properties. If $v$ is a concave vertex and $w$
is not adjacent to $v$, and yet is visible from $v$, then
$|v-w|\ge3.2$. If $v$ is a concave corner, then
    $|v-w|\ge3.07$ for both adjacent corners $w$.
If $v$ is a concave corner and $|v|\ge2.2$, then
    $|v-w|\ge3.2$ for both adjacent corners $w$.
These hypotheses will remain in force through the end of
Section \refz 4.

We call a spherical region convex if its interior angles are all
less than $\pi$. The case where the subregion is a convex triangle
will be treated in Section \refz 5.7. Hence, we may also assume in
Sections \refz 4.10 through \refz 4.13 that the subregion is not a
convex triangle.

We construct a {\it corner cell\/}
at each corner.  It depends
on a parameter $\lambda \in [1.6,1.945]$.
In all applications, we take $\lambda = 1.945 = 3.2-t_0$,
$\lambda = 1.815 = 3.07-t_0$, or $\lambda = 1.6 = 3.2/2$.

To construct the cell around
the corner $v$,
place a triangle along $(0,v)$ with sides $|v|$, $t_0$, $\lambda$
(with $\lambda$ opposite the origin).
Generate the solid of rotation around the axis $(0,v)$.  Extend to
a cone over $0$.  Slice the solid by the perpendicular bisector of
$(0,v)$, retaining the part near $0$.  Intersect the solid with
a ball of radius $t_0$.   The cones over the two boundary edges
of the subregion at $v$ make two cuts in the solid.  Remove the
slice that lies outside the cone over the subcluster.  What remains
is the corner cell at $v$ with parameter $\lambda$.

Corner cells at corners separated
by a distance less than $2\lambda$ may overlap.  We define a truncation of
the corner cell that has the property that the {\it
truncated corner cells\/}
at adjacent corners do not overlap.
Let $(0,v_i,v_j)^\perp$ denote the plane perpendicular to the
plane $(0,v_i,v_j)$ passing through the origin and the
circumcenter of $(0,v_i,v_j)$.

Let $v_1,v_2,v_3$ be consecutive corners of a subcluster. Take the
corner cell with parameter $\lambda$ at the corner $v_2$.  Slice
it by the planes $(0,v_1,v_2)^\perp$ and $(0,v_2,v_3)^\perp$, and
retain the part along the edge $(0,v_2)$. This is the truncated
corner cell (tcc).  By construction tccs at adjacent corners are
separated by a plane $(0,\cdot,\cdot)^\perp$. Tccs at nonadjacent
corners do not overlap if the corners are $\ge2\lambda$ apart.
Tccs will only be used in subregions satisfying this condition. It
will be shown in Section \refz 4.12 that tccs lie in the cone over
the subregion (for suitable $\lambda$).

\subhead \refz 4.11 Formulas for Truncated corner cells\endsubhead

We will assign a score to truncated corner cells, in such a way
that the score of the subcluster can be estimated from the scores
of the corner cells.

We write $C_0$ for a truncated corner cell.  We write $C_0^u$ for
the corresponding untruncated corner cell.  (Although we call this
the untruncated corner cell to distinguish it from the corner cell,
it is still truncated in the sense that it lies in the ball at the origin
of radius $t_0$.  It is untruncated in the sense that it is not
cut by the planes $(\ldots)^\perp$.)

For any solid body $X$, we define the {\it geometric}
truncated Voronoi function
by
    $$\vor_0^{g}(X) = 4(-\doct \Vol(X) + \sol(X)/3)$$
the counterpart for squander
    $$\tau_0^g(X) = \zeta\pt \sol(X) - \vor_0^g(X).$$
The solid angle is to be interpreted as the solid angle of the
cone formed by all rays from the origin through nonzero points of
$X$. We may apply these definitions to obtain formulas for
$\vor_0^{g}(C_0)$, and so forth.

The formula for the score of a truncated corner cell differs slightly
according to the convexity of the corner.  We start with a convex
corner $v$, and let $v_1$, $v$, and $v_2$ be consecutive corners
in the subregion.

Let $S=(0,v,v_1,v_2)$ be a simplex with $|v_1-v_2|\ge3.2$.
The formula for the score of a tcc $C_0(S)$ simplifies if the
face of $C_0$
cut by $(0,v,v_1)^\perp$ does not meet the
face cut by $(0,v,v_2)^\perp$.
We make that assumption in this subsection.
  Set $\chi_0(S) = \vor^g_0(C_0(S))$.
(The function $\chi_0$ is unrelated to the function $\chi$ that was
introduced in Section I.8 to measure the orientation of faces.)

$$\align
\psi &= \arc(y_1,t_0,\lambda),\quad h=y_1/2,\\
R'_{126}&=R(y_1/2,\eta_{126},y_1/(2\cos\psi)),
\quad R_{126}=R(y_1/2,\eta_{126},t_0),\\
\sol'(y_1,y_2,y_6) &= +\dih(R'_{126})(1-\cos\psi)-\sol(R'_{126}),\\
\chi_0(S) &= \dih(S)(1-\cos\psi)\phi_0\\
    &\ -\sol'(y_1,y_2,y_6)\phi_0 -\sol'(y_1,y_3,y_5)\phi_0\\
    &\ +A(h)\dih(S)-
        4\doct (\quo(R_{126})+\quo(R_{135})).
\endalign
$$
In the three lines giving the formula for $\chi_0$, the first line
    represents the score of the cone before it is cut by the
    planes $(0,v,v_i)^\perp$ and the perpendicular bisector of $(0,v)$.
The second line is the correction resulting from cutting the tcc
    along the planes $(0,v,v_i)^\perp$.
    The face of the Rogers simplex $R'_{126}$ lies along the
    plane $(0,v,v_1)^\perp$.  The third line is the
    correction from slicing the tcc with the perpendicular bisector
    of $(0,v)$.  This last term is the same as the term appearing
    for a similar reason in the formula for $\vor_0$ in F.3.7.
    In this formula $R$ is the usual Rogers simplex and $\quo(R_{ijk})$
    is the quoin coming from a Rogers simplex along the face with
    edges $(ijk)$.

The formula for the untruncated corner cell is obtained by setting
``$\sol'$'' and ``$\quo$'' to ``$0$'' in the expression for $\chi_0$.
Thus,
$$
\vor^g(C_0^u) = \dih(S)[(1-\cos\psi)\phi_0 + A(h)]
$$
The formula depends only on $\lambda$, the dihedral angle, and the
height $|v|$.  We write $C_0^u = C_0^u(|v|,\dih)$, and suppress
$\lambda$ from the notation. The dependence on $\dih(S)$ is
linear:
 $$
\tau^g_0(C_0^u(|v|,\dih))= (\dih/\pi)\tau^g_0(C_0^u(|v|,\pi)).
$$

The dependence of $\chi_0$ on the fourth edge $y_4=|v_1-v_2|$
comes through a term proportional to $\dih(S)$.  Since the
dihedral angle is monotonic in $y_4$, so is $\chi_0$.  Thus, under
the assumption that $|v_1-v_2|\ge3.2$,  we obtain an upper bound
on $\chi_0$ at $y_4=3.2$. Our deformations will fix the lengths of
the other five variables, and monotonicity gives us the sixth.
Thus, the tccs lead to an upper bound on $\vor^g_0$ (and a lower
bound on $\tau^g_0$) that does not require interval arithmetic.

At a concave vertex, the formula is similar.  Replace ``$\dih(S)$''
with $``(2\pi-\dih(S))$'' in the given expression for $\chi_0$.
We add a superscript $-$ to the name of the function at concave
vertices, to denote this modification: $\chi^-_0(C_0)$.

\subhead \refz 4.12 Containment of Truncated corner
cells\endsubhead
\smallskip

The assumptions made at the beginning of Section \refz 4.10 remain in
force.

\proclaim{Lemma \refz 4.12.1} Let $v$ be a concave vertex with
$|v|\ge2.2$. The truncated corner cell at $v$ with parameter
$\lambda=1.945$ lies in the truncated $V$-cell over $R$.
\endproclaim

\demo{Proof}
Consider a  corner cell at $v$ and a distinguished
edge $(v_1,v_2)$ forming the boundary of the subregion.
The corner cell with parameter $\lambda=1.945$
is contained in a cone of arcradius
$\theta = \arc(2,t_0,\lambda)< 1.21 <\pi/2$
(in terms of the function {\it arc\/} of Section 2.8).
Take two corners $w_1$, $w_2$, visible from $v$,
 between which the given bounding edge
appears. (We may have $w_i=v_i$). The two visible edges,
$(v,w_i)$,  have length $\ge 3.2$. (Recall that the distinguished
edges at $v$ have been deformed to length $3.2$.) They have
arc-length at least $\arc(2.51,2.51,3.2)>1.38$. The segment of the
distinguished edge $(v_1,v_2)$ visible from $v$ has arc-length at
most $\arc(2,2,3.2)<1.86$.

We check that the corner cell cannot cross the visible
portion of the edge $(v_1,v_2)$.
Consider the spherical triangle formed by the edges $(v,w_1)$,
$(v,w_2)$ (extended as needed) and the visible part of $(v_1,v_2)$.
Let $C$ be the projection of
$v$ and $AB$ be the projection of the visible part of $(v_1,v_2)$.
Pivot $A$ and $B$ toward $C$ until the edges $AC$ and $BC$ have
arc-length $1.38$.  The perpendicular
from $C$ to $AB$ has length at least
    $$\arccos(\cos(1.38)/\cos(1.86/2))>1.21>\theta.$$
This proves that the corner cell lies in the cone over the subregion.
\qed\enddemo

\proclaim{Lemma \refz 4.12.2} Let $v$ be a concave vertex. The
truncated corner cell at $v$ with parameter $\lambda=1.815$ lies
in the truncated $V$-cell over $R$.
\endproclaim

\demo{Proof}
The proof proceeds along the same lines as the previous lemma, with
slightly different constants.
Replace $1.945$ with $1.815$, $1.38$ with $1.316$, $1.21$ with $1.1$.
Replace $3.2$ with $3.07$ in contexts giving a lower bound to the
length of an edge at $v$, and keep it at $3.2$ in contexts calling
for an upper bound on the length of a distinguished edge.
The constant $1.86$ remains unchanged.
\qed
\enddemo

\proclaim{\bf Lemma \refz 4.12.3} The truncated corner cells with
parameter $1.6$ in a subregion do not overlap.
\endproclaim

\demo{Proof}
We may assume that the corners are not adjacent.
If a nonadjacent corner $w$ is visible from $v$, then $|w-v|\ge3.2$,
and an interior point intersection $p$ is incompatible with the triangle
inequality: $|p-v|\le 1.6$, $|p-w|<1.6$.
If $w$ is not visible, we have a chain $v=v_0,v_1,\ldots,v_r=w$ such
that $v_{i+1}$ is visible from $v_i$. Imagine a taut string
inside the subregion extending from $v$ to $w$.  The projections
of $v_i$ are the
corners of the string's path.   The string bends in an angle greater
than $\pi$ at each $v_i$, so the angle at each intermediate
$v_i$ is greater than
$\pi$. That is,
they are concave. Thus, by our deformations $|v_i-v_{i+1}|\ge3.07$.
  The
string has arc-length at least $r \arc(2.51,2.51,3.07)>r (1.316)$.
But the corner cells
lie in cones of arcradius $\arc(2,t_0,\lambda)< 1$.
So $2(1.0)>r(1.316)$, or $r=1$.  Thus, $w$ is visible from $v$.
\qed
\enddemo

\proclaim{\bf Lemma \refz 4.12.4} The corner cell for $\lambda \le
1.815$ does not overlap the $t_0$-cone wedge around another corner
$w$.
\endproclaim

\demo{Proof}
We take $\lambda=1.815$.
As in the previous proof, if there is overlap along a chain,
then
$$\arc(2,t_0,\lambda) +\arc(2,t_0,t_0) > r \arc(2.51,2.51,3.07),$$
and again $r=1$.  So each of the two vertices in question is
visible from the other. But overlap implies $|p-v|\le1.815$ and
$|p-w|<1.255$, forcing the contradiction $|w-v|<3.07$. \qed
\enddemo

\proclaim{\bf Lemma \refz 4.12.5} The corner cell for $\lambda \le
1.945$ at a corner $v$ satisfying $|v|\ge2.2$ does not overlap the
$t_0$-cone wedge around another corner $w$.
\endproclaim

\demo{Proof}
We take $\lambda=1.945$.
As in the previous proof, if there is overlap along a chain, then
$$\arc(2,t_0,\lambda) +\arc(2,t_0,t_0) > r \arc(2.51,2.51,3.2),$$
and again $r=1$.  Then the result follows from
$$|w-v|\le |p-v|+|p-w| < 1.945 + 1.255 = 3.2.$$
\qed
\enddemo

Lemma \refz 4.3 was stated in the context of a subregion before
deformation, but a cursory inspection of the proof shows that the
geometric conditions required for the proof remain valid by our
deformations. (This assumes that the subregion is not a triangle,
which we assumed at the beginning of Section \refz 4.10.) In more
detail, there is a solid $CP_0$ contained in the ball of radius of
$t_0$ at the origin, and lying over the cone of the subregion $P$
such that a bound on the penalty-free subcluster score is
$\vor^g_0(CP_0)$ and squander $\tau^g_0(CP_0)$. (By {\it
penalty-free\/} score, we mean the part of the scoring bound that
does not include any of the penalty terms.  We will sometimes call
the full score, including the penalty terms, the {\it
penalty-inclusive\/} score.)

Let $\{y_1,\ldots,y_r\}$ be a decomposition of the subregion into disjoint
regions whose union is $X$.
Then if we let $CP_0(y_i)$ denote the intersection of $CP_0(y_i)$ with
the cone over $y_i$, we can write
    $$\tau^g_0(CP_0) =\sum_i \tau^g_0(CP_0(y_i)).$$

These lemmas allow us to express bounds on the score (and
squander) of a subcluster as a sum of terms associated with
individual (truncated) corner cells. By Lemmas \refz 4.12.1
through \refz 4.12.5, these objects do not overlap under suitable
conditions. Moreover, by the interpretation of terms provided by
Section \refz 4.3, the cones over these objects do not overlap,
when the objects themselves do not. In other words, under the
various conditions, we can take the (truncated) corner cells to be
among the sets $CP_0(y_i)$.

To work a typical example, let us place a truncated corner cell
with parameter $\lambda=1.6$ at each concave corner.  Place a
$t_0$-cone wedge $X_0$ at each convex corner. The cone over each
object lies in the cone over the subregion. By Lemma \refz 4.3 and
Lemma \refz 4.1 (see the proof), the $t_0$-cone wedge $X_0$
squanders a positive amount.  The part $P'$ of the subregion
outside all truncated corner cells and outside the $t_0$-cone
wedges squanders
$$\sol(P')(\zeta\pt-\phi_0) > 0.$$
where $\sol(P')$ is the part of the solid angle of the subregion
lying outside the tccs. Dropping these positive terms, we obtain a
lower bound on the penalty-free squander:
    $$\tau^g_0(CP_0) \ge \sum_{C_0} \tau^g_0(C_0).$$
There is one summand for each concave corner of the subregion.
Other cases proceed similarly.

\subhead \refz 4.13 Convexity\endsubhead

\proclaim{Lemma 4.13.1} There are at most two concave corners.
\endproclaim

\demo{Proof}
Use the parameter $\lambda=1.6$ and place a truncated corner cell
$C_0$ at
each concave corner $v$.
Let $C_0^u(|v|,\dih)$ denote the corresponding untruncated cell.
Formula \refz 4.11  gives
    $$
    \tau_0^g(C_0) =
    \tau^g_0(C^u_0(|v|,\dih))
    - \sol'(y_1,y_2,y_6) \phi'_0
    -\sol'(y_1,y_3,y_5)\phi'_0,
    $$
where $\phi'_0 = \zeta\pt-\phi_0 < 0.6671$.
(The conditions $y_5\ge3.07$ and $y_6\ge 3.07$
force the faces along the these edges to have circumradius greater than
$t_0$, and this causes the ``$\quo$'' terms in the formula to be zero.)

By monotonicity in $\dih$,
a lower bound on $\tau^g_0(C^u_0)$ is obtained at $\dih=\pi$.
$\tau_0(C^u_0(|v|,\pi))$ is an explicit monotone decreasing
rational function of $|v|\in[2,2.51]$,
which is minimized for $|v|=2.51$.  We find
$$\tau_0(C_0^u(|v|,\dih))\ge\tau_0(C_0^u(2.51,\pi)) >0.32.$$

The term $\sol'(y_1,y_3,y_5)$ is maximized when $y_3=2.51$,
$y_5=3.07$, so that $\sol'< 0.017$.  (This was checked with
interval arithmetic in Mathematica.) Thus,
    $$\tau_0(C_0(v))\ge 0.32 - 2(0.017) \phi_0' > 0.297.$$

If there are three or more concave corners,
then the penalty-free corner cells squander
at least $3(0.297)$.
The penalty is at most $\maxpi$ (Section 4.7).
So the penalty-inclusive squander
is more than $3(0.297) - \maxpi >\squander$.
\qed
\enddemo

\proclaim{Lemma 4.13.2}
    There are no concave corners of height at most $2.2$.
\endproclaim

\demo{Proof} Suppose there is a corner of height at most $2.2$.
Place an untruncated corner cell $C^u_0(|v|,\dih)$ with parameter
$\lambda =1.815$ at that corner and a $t_0$-cone wedge at every
other corner.  The subcluster squanders at least
    $\tau_0(C_0(|v|,\pi))-\maxpi$.
This is an explicit monotone decreasing rational function of one
variable. The penalty-inclusive squander is at least
    $$\tau_0(C^u_0(2.51,\pi))-\maxpi > \squander.$$
\qed
\enddemo

By the assumptions at the beginning of Section \refz 4.10,
the lemma implies that
each concave corner has distance at least $3.2$ from every
other visible corner.

As in the previous lemma, when $\lambda=1.945$, a lower bound on
what is squandered by the corner cell is obtained for  $|v|=2.51$,
$\dih=\pi$.  The explicit formulas give
penalty-free squander $>0.734$.
Two disjoint corner cells give
penalty-inclusive squander $>\squander$.  Suppose two at
$v_1,v_2$ overlap.  The lowest bound is obtained when $|v_1-v_2|=3.2$,
the shortest distance possible.

We define a function $f(y_1,y_2)$
that measures what the union of the overlapping
corner cells squander.  Set $y_i = |v_i|$, $\ell=3.2$, and
$$
\align
\alpha_1 &= \dih(y_1,t_0,y_2,\lambda,\ell,\lambda),\\
\alpha_2 &= \dih(y_2,t_0,y_1,\lambda,\ell,\lambda),\\
\sol &= \sol(y_2,t_0,y_1,\lambda,\ell,\lambda),\\
\phi_i &= \phi(y_i/2,t_0),\quad i=1,2,\\
\lambda&=3.2-t_0=1.945,\\
f(y_1,y_2)&=
    2(\zeta\pt-\phi_0)\sol+
    2\sum_1^2 \alpha_i(1-y_i/(2t_0))(\phi_0-\phi_i)\\
        &\quad +\sum_1^2 \tau_0(C(y_i,\lambda,\pi-2\alpha_i)).
\endalign
$$
$\A_{14}$ gives $f(y_1,y_2)>\squander+\maxpi$,
for $y_1,y_2\in[2,2.51]$.

We conclude that there is at most one concave corner. Let $v$ be
such a corner.  If we push $v$ toward the origin, the solid angle
is unchanged and $\vor_0$ is increased.  Following this by the
deformation of Section 4.9, we maintain the constraints
$|v-w|=3.2$, for adjacent corners $w$, while moving $v$ toward the
origin. Eventually $|v|=2.2$.  This is impossible by Lemma \refz
4.13.2.

We verify that this deformation preserves the constraint $|v-w|\ge2$,
for all corners $w$ such that $(v,w)$ lies entirely outside the
subregion.  If fact,  every corner
is visible from $v$, so that the subregion is star convex at
$v$.
  We leave the details to the reader.

We conclude that all subregions
can be deformed into convex polygons.

\heads{\refz 5. Convex Polygons}
\vfill\eject
\head \refz 5. Convex Polygons\endhead

\subhead\refz 5.1 Deformations\endsubhead We divide the bounding
edges over the polygon according to length $[2,2.51]$,
$[2.51,2\sqrt{2}]$, $[2\sqrt{2},3.2]$. The deformations of Section
4.9 contract edges to the lower bound of the intervals ($2,2.51$,
or $2\sqrt{2}$) unless a new distinguished edge is formed.  By
deforming the polygon, we assume that the bounding edges have
length $2,2.51$, or $2\sqrt{2}$. (There are a few instances of
triangles or quadrilaterals that do not satisfy the hypotheses
needed for the deformations. These instances will be treated in
Sections 5.7 and 5.8.)

\proclaim{Lemma \refz 5.1.1}  Let $S=S(y_1,\ldots,y_6)$ be a simplex,
with $x_i=y_i^2$, as usual.  Let $y_4\ge 2$, $\Delta\ge0$,
$y_5,y_6\in\{2,2.51,2\sqrt{2}\}$.  Fixing all the variables but $x_1$,
let $f(x_1)$ be one of
the functions $\vor_0(S)$ or $-\tau_0(S)$. We have
$f''(x_1)>0$ whenever $f'(x_1)=0$.
\endproclaim

\demo{Proof} This is an interval calculation $\A_{15}$.\qed\enddemo

The lemma implies that $f$ does not have an interior point local
maximum for $x_1\in[2^2,2.51^2]$.  Fix three consecutive corners,
$v_0,v_1,v_2$  of the convex polygon, and apply the lemma to the
variable $x_1 = |v_1|^2$ of the simplex $S=(0,v_0,v_1,v_2)$. We
deform the simplex, increasing $f$.  If the deformation produces
$\Delta(S)=0$, then some dihedral angle is $\pi$, and the
arguments for nonconvex regions bring us eventually back to the
convex situation. Eventually $y_1$ is $2$ or $2.51$.  Applying the
lemma at each corner, we may assume that the height of every
corner is $2$ or $2.51$.   (There are a few cases where the
hypotheses of the lemma are not met, and these are discussed in
Sections 5.7 and 5.8.)

\proclaim{Lemma \refz 5.1.2} The convex polygon has at most $7$ sides.
\endproclaim

\demo{Proof}  Since the polygon is convex, its perimeter on the
unit sphere is
at most a great circle $2\pi$.  If there are $8$ sides, the
perimeter is at least $8\arc(2.51,2.51,2)>2\pi$.
\qed
\enddemo

\subhead \refz 5.2 Truncated corner cells\endsubhead

The following lemma justifies using tccs at the corners as an
upper bound on the score (and lower bound on what is squandered).
We fix the truncation parameter at $\lambda=1.6$.

\proclaim{Lemma}
Take a convex subregion that is not a triangle.  Assume
edges between adjacent corners have lengths $\in\{2,2.51,2\sqrt{2},3.2\}$.
Assume nonadjacent
corners are separated by distances $\ge3.2$.
 Then the truncated corner cell at each vertex lies in the
cone over the subregion.
\endproclaim

\demo{Proof} Place a tcc at $v_1$. For a contradiction,
let $(v_2,v_3)$ be an edge that the tcc overlaps.  Assume first
that $|v_1-v_i|\ge 2.51$, $i=2,3$.  Pivot so that
$|v_1-v_2|=|v_1-v_3|=2.51$.    Write
$S(y_1,\ldots,y_6)=(0,v_1,v_2,v_3)$.
Set $\psi=\arc(y_1,t_0,1.6)$. $\A_1$ gives
$\beta_\psi(y_1,y_2,y_6)<\dih_2(S)$.

Now assume $|v_1-v_2|<2.51$.  By the hypotheses of the lemma,
$|v_1-v_2|=2$.  If $|v_1-v_3|<3.2$, then $(0,v_1,v_2,v_3)$
is triangular, contrary to hypothesis.  So $|v_1-v_3|\ge3.2$.
Pivot so that  $|v_1-v_3|=3.2$. By $\A_1$,
    $$\beta_\psi(y_1,y_2,y_6)< \dih_2(S),$$
where $\psi=\arc(y_1,t_0,1.6)$, provided $y_1\in[2.2,2.51]$.
Also,
if $y_1\in[2.2,2.51]$
$$\arc(y_1,t_0,1.6)<\arc(y_1,y_2,y_6).$$

If $y_1\le 2.2$, then $\Delta_1\ge0$, so
$\partial\dih_2/\partial x_3\le 0$.  Set $x_3=2.51^2$.
Also, $\Delta_6\ge0$, so $\partial\dih_2/\partial x_4\le0$.
Set $x_4=3.2^2$.

Let $c$ be a point of intersection of the plane $(0,v_1,v_2)^\perp$
with the circle at distance $\lambda=1.6$ from $v_1$ on the sphere
centered at the origin
of radius $t_0$.  The angle along $(0,v_2)$ between the planes
$(0,v_2,v_1)$ and $(0,v_2,c)$ is
    $$\dih(R(y_2/2,\eta_{126},y_1/(2\cos\psi))).$$
This angle is less than $\dih_2(S)$ by $\A_1$.
Also, $\Delta_1\ge0$, $\partial\dih_3/\partial x_2\le0$,
so set $x_2=2.51^2$.
Then $\Delta_5<0$, so $\dih_2>\pi/2$.  This means that $(0,v_1,v_2)^\perp$
separates the tcc from the edge $(v_2,v_3)$.
\qed
\enddemo

\subhead \refz 5.3 Analytic continuation\endsubhead

In this subsection we assume that $\lambda=1.6$ and that the
truncated corner cell under consideration lies at a convex vertex.

Assume that the face cut by $(0,v,v_1)^\perp$ meets the face
cut by $(0,v,v_2)^\perp$.  Let $c_i$ be the point on the
plane $(0,v,v_i)^\perp$ satisfying $|c_i-v|=1.6$, $|c_i|=t_0$.
(Pick the root within the wedge between $v_1$ and $v_2$.)
The overlap of the two faces is represented in the diagram.

\smallskip
\gram|1.8|5.3|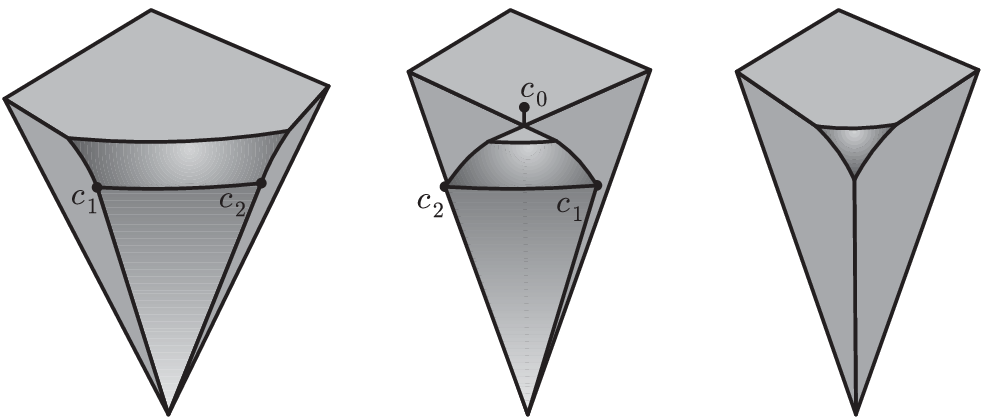|
\smallskip

We let $c_0$ be the point of height $t_0=1.255$
on the intersection of the planes $(0,v,v_1)^\perp$ and $(0,v,v_2)^\perp$.
We claim
that $c_0$ lies over the truncated spherical region
of the tcc, rather than the wedges of $t_0$-cones or the Rogers simplices
along the faces $(0,v,v_1)$ and $(0,v,v_2)$.  (This implies that $c_0$
cannot protrude beyond the corner cell as depicted in the second frame
of the diagram.)
To see the claim, consider
the tcc as a function of $y_4=|v_1-v_2|$.  When $y_4$ is sufficiently
large the claim is certainly true.  Contract $y_4$ until $c_0=c_0(y_4)$ meets
the perpendicular bisector of $(0,v)$.  Then $c_0$ is equidistant
from $0,v,v_1$ and $v_2$ so it is the circumcenter of $(0,v,v_1,v_2)$.
It has distance $t_0$ from the origin, so the circumradius is $t_0$.
This implies that $y_4\le 2.51$.

The tcc is defined by the constraints represented
in the third frame.  The analytic continuation of the function
$\chi_0(S)=\chi_0^\anal(S)$, defined above,
acquires a volume $X$, counted
with negative sign, lying under the spherical triangle $(c_0,c_1,c_2)$.
  Extending our
notation, we have an analytically defined function $\chi_0^\anal$
 and
a geometrically defined function $\chi_0^\geom$,
$$\align
\chi_0^\anal(S) &= \chi_0^\geom(S)-\vor_0(X), \text{\ where}\\
\vor_0(X) &= 4(-\doct\Vol(X) +\sol(X)/3) = \phi_0\sol(X) <0.
\endalign
$$
So $\chi_0^\anal >\chi_0^\geom$, and we may always use
$\chi_0(S)=\chi_0^\anal(S)$ as an upper bound on the score of a tcc.

For example, with $\lambda=1.6$
and $S = S(2.3,2.3,2.3,2.9,2,2)$, we have
$$\chi_0^\anal(S)\approx -0.103981, \quad\chi_0^\geom(S)\approx -0.105102.$$
Or, if $S=S(2,2,2.51,3.2,2,2.51)$, then
$$\chi_0^\anal(S)\approx -0.0718957, \quad
\chi_0^\geom(S)\approx -0.0726143.$$

\subhead \refz 5.4 Penalties\endsubhead

In Section 4.7, we determined the bound of $\maxpi=0.06688$ on
penalties. In this section, we give a more thorough treatment of
penalties. Until now a penalty has been associated with a given
standard region, but by taking the worst case on each subregion,
we can move the penalties to the level of subregions.   Roughly,
each subregion should incur the penalties from the upright
quarters that were erased along edges of that subregion.  Each
upright quarter of the original standard region is attached at an
edge between adjacent corners of the standard cluster. The edges
have lengths between $2$ and $2.51$.  The deformations shrink the
edges to length $2$.  We attach the penalty from the upright
quarter to this edge of this subregion. In general, we divide the
penalty evenly among the upright quarters along a common diagonal,
without trying to determine a more detailed accounting. For
example, the penalty $0.008$ in Section 3.9 comes from three
upright quarters.  Thus, we give each of three edges a penalty of
$0.008/3$.  Or, if there are only two upright quarters in the
group $\Splus$, then each of the two upright quarters is assigned
the penalty $0.00222/2$ (see Lemma \refz 3.9.2).

The penalty $0.04683 = 3\xiG$ in Section 4.7 comes from
three upright quarters $\Sminus$.   Each of
three edges is assigned a penalty of $\xiG$.  The penalty
$0.03344=3\xiG+\xikG$
comes from the arrangement of four upright quarters
$\Sfour$ of Section 3.8.
It is divided among 4 edges. These are the
only upright quarters that take a penalty when erased.
(The case of two upright quarters over a flat quarter as in Lemma 3.4,
are treated by a separate argument in Section 5.7.
Loops will be discussed in Section 5.11.)

The penalty can be reduced in various situations involving a
masked flat quarter.  For example, in the three-quarter
configuration $\Sminus$, if there is a masked flat quarter, two of
the uprights are scored by the analytic Voronoi function, so that
the penalty plus adjustment is only $0.034052=2\xiV+\xiG+0.0114$
(by $\A_{10},\A_{11}$).  The adjustment $0.0114$ reflects the
scoring rules for masked flat quarters (Section 3.9).  This we
divide evenly among the three edges that carried the upright
quarters. If $e$ is an edge of the subregion $R$, let $\pi_0(R,e)$
denote the penalty and score adjustment along edge $e$ of $R$.

In summary,
we have the penalties,
    $$\xik,\xiV,\xiG,\ 0.008,$$
combined in various ways in the configurations $\Sminus$, $\Splus$,
$\Sfour$.  There are score adjustments
$$0.0114\quad \text{ and }\quad 0.0063$$
from Section 3.10 for masked flat quarters.  If the sum of
these contributions is $s$, we set $\pi_0(R,e)=s/n$, for
each edge $e$ of $R$ originating from an erased upright quarter
of $\Cal{\bold S}_n^\pm$.

\subhead \refz 5.5 Penalties and Bounds\endsubhead

Recall that the bounds for flat quarters we wish to establish from
Section 4.5 are
$Z(3,1)=0.00005$ and $D(3,1)=0.06585$.
Flat quarters arise
in two different ways.  Some flat quarters are present before the
deformations begin.  They are scored by the rules of Section 3.10.
Others are formed by the deformations.  In this case, they are
scored by $\vor_0$.
Since the flat quarter is broken away from the subregion
 as soon as the diagonal reaches
$2\sqrt{2}$, and then is not deformed further, the diagonal
is fixed at $2\sqrt{2}$.  Such flat quarters can violate our
desired inequalities. For example,
$$Z(3,1)<\vor_0(S(2,2,2,2\sqrt{2},2,2)) \approx 0.00898,\quad
    \tau_0(S(2,2,2,2\sqrt{2},2,2))\approx 0.0593.$$
On the other hand, as we will see, the adjacent subregion satisfies the
inequality by a comfortable margin.  Therefore, we define a transfer
$\epsilon$ from flat quarters to the adjacent subregion.
(In an exceptional region,
the subregion next to a flat quarter along the diagonal
is not a flat quarter.)

For a flat quarter $Q$, set
$$
\epsilon_\tau(Q) = \cases 0.0066,&\text{(deformation),}\\
            0,&\text{(otherwise)}.
            \endcases
$$
$$
\epsilon_\sigma(Q) = \cases 0.009,&\text{(deformation),}\\
            0,&\text{(otherwise)}.
            \endcases
$$
The nonzero value occurs when the flat quarter $Q$ is obtained by
deformation from an initial configuration in which $Q$ is not
a quarter.
The value is zero when the flat quarter $Q$ appears already
in the undeformed standard cluster.
Set
$$\align
\pi_\tau(R) &= \sum_e \pi_0(R,e) +
    \sum_e\pi_0(Q,e)+\sum_Q \epsilon_\tau(Q),\\
\pi_\sigma(R)&=\sum_e \pi_0(R,e) +
    \sum_e\pi_0(Q,e)+\sum_Q \epsilon_\sigma(Q).
\endalign
$$
The first sum runs over the edges of
a subregion $R$.  The second sum
runs over the edges of the
flat quarters $Q$ that lie adjacent to $R$ along the diagonal
of $Q$.

The edges between corners of the polygon have lengths $2$, $2.51$,
or $2\sqrt{2}$.  Let $k_0$, $k_1$, and $k_2$ be the number
of edges of these three lengths respectively.  By Lemma 5.1, we
have $k_0+k_1+k_2\le7$.  Let $\tilde\sigma$ denote any of the functions
of Section 3.10.(a)--(f).  Let $\tilde\tau = \sol\zeta\pt - \tilde\sigma$.

To prove Theorem 4.4, refining the strategy proposed in Section 4.5,
we must show that for each flat quarter $Q$ and each
subregion $R$ that
is not a flat quarter, we have
$$
\align
\tilde\tau(Q) &> D(3,1) - \epsilon_\tau(Q),\\
\tau_0(Q) &> D(3,1)-\epsilon_\tau(Q),\quad\text{if }y_4(Q)=2\sqrt2,\\
\tau_V(R) &> D(3,2),\quad\text{(type $S_A$)},\\
\tau_0(R) &> D(k_0+k_1+k_2,k_1+k_2)+\pi_\tau(R),\tag 5.5.1
\endalign$$
where $D(n,k)$ is the function defined in Section 4.5.
The first of these inequalities follows from
$\A_{1},\A_{13},\A_{16}$.
In general, we are given a subregion without explicit information
about what the adjacent subregions are.  Similarly, we have
discarded all information about what upright quarters have been
erased.  Because of this, we
assume the worst, and use the largest feasible values of $\pi_\tau$.

\proclaim{Lemma}  We have
$\pi_\tau(R)\le 0.04683 + (k_0+2k_2-3)0.008/3 +0.0066k_2$.
\endproclaim

\demo{Proof}
The worst penalty $0.04683=3\xiG$ per edge comes from $\Sminus$.
The number of penalized edges not on $\Sminus$ is at most
$k_0+2k_2-3$.
For every three edges we might have one $\Splus$.
The other cases such as $\Sfour$ or situations
with a masked flat quarter are readily seen to give smaller penalties.
\qed\enddemo

For bounds on the score, the situation is similar.  The
only penalties we need to consider are $0.008$ from Section 3.9.
If either of the other configurations of upright quarters $\Sminus$,
$\Sfour$ occur, then the
score of the standard cluster is less than $s_8=-0.228$,
by Sections 3.7 and 3.8. This is the desired bound.  So it is enough to
consider subregions that do not have these upright configurations.  Moreover,
the penalty $0.008$ does not occur in connection with masked
flats.  So we can take $\pi_\sigma(R)$ to be
    $$(k_0+2k_2)0.008/3 + 0.009 k_2.$$
If $k_0+2k_2<3$, we can strengthen this to
$\pi_\sigma(R)=0.009 k_2$.
Let $\tilde\sigma$ be any of the functions of Section 3.10.(a)--(f).
To prove Theorem 4.4, we will show
    $$
\align
\tilde\sigma(Q)&< Z(3,1)+\epsilon_\sigma(Q),\\
\vor_0(Q)&< Z(3,1)+\epsilon_\sigma(Q),\quad\text{if }y_4(Q)=2\sqrt2,\\
\vor_0(R)&< Z(3,2),\quad\text{(type $S_A$)},\\
\vor_0(R) &< Z(k_0+k_1+k_2,k_1+k_2) - \pi_\sigma(R).
    \tag 5.5.2
\endalign
$$
The first of these inequalities follows form $\A_1,\A_{13},\A_{16}$.

\subhead \refz 5.6 Constants\endsubhead

Theorem 4.4 now results from the calculation of a host of
constants.  Perhaps there are simpler ways to do it, but it was a
routine matter to run through the long list of constants by
computer. What must be checked is that the Inequalities \refz
5.5.1 and \refz 5.5.2 hold for all possible convex subregions.
This section describes in detail the constants to check.

We begin with a subregion given as a convex $n$-gon, with at least
4 sides.   The heights
of the corners and the lengths of edges between adjacent edges
have been reduced by deformation to a finite number of possibilities
(lengths 2,2.51, or lengths 2,2.51,$2\sqrt{2}$, respectively).
By Lemma 5.1, we may take $n=4,5,6,7$.   Not all possible assignments of
lengths correspond to a geometrically viable configuration.
One constraint that eliminates many possibilities, especially heptagons,
is that of Section 5.1: the perimeter of the convex
polygon is at most
a great circle.  Eliminate all length-combinations that do not
satisfy this condition.  When there is a special simplex it can be
broken from the subregion and scored separately unless the two heights
along the diagonal are 2 (see $\A_{13}$).
We assume in all that follows that
all specials that can be broken off have been.
There is a second condition related to
special simplices.  We have $\Delta(2.51^2,2^2,2^2,x^2,2^2,2^2)<0$,
if $x> 3.114467$.  This means that if the cluster edges
along the polygon are $(y_1,y_2,y_3,y_5,y_6)= (2.51,2,2,2,2)$,
the simplex must be special ($y_4\in[2\sqrt{2},3.2]$).

The easiest cases
to check are those with no special simplices over the polygon.  In other
words, these are subregions for which the distances between nonadjacent
corners are at least 3.2.  In this case we approximate the score
(and what is squandered) by tccs at the corners.  We use
monotonicity to bring the fourth edge to length $3.2$.
We calculate the tcc constant bounding the score, checking that
it is less than the constant
    $ Z(k_0+k_1+k_2,k_1+k_2) - \pi_\sigma$, from
(5.5.2).
The bounds for $\tau_0$ are verified in the same way.

When $n=5,6,7,$ and there is one special simplex, the situation is not
much more difficult.  By our deformations,  we decrease
the lengths of edges $2,3,5,6$ of the special to 2.
We remove the special by cutting along
 its fourth edge $e$ (the diagonal).
  We score the special with the weak bounds found in
$\A_{13}$.  Along the edge $e$, we then apply deformations to the
$(n-1)$-gon that remains.
If this deformation brings $e$ to
length $2\sqrt{2}$, then the $(n-1)$-gon may be scored
with tccs as in the previous paragraph.  But there are other possibilities.
Before $e$ drops to $2\sqrt{2}$, a new distinguished edge of length
$3.2$ may form between two corners (one of the corners will be a
chosen endpoint of $e$).  The subregion breaks in two.
By deformations, we eventually arrive at
 $e=2\sqrt2$ and a subregion with diagonals
of length at least $3.2$.  (There is one case that
may fail to be deformable to $e=2\sqrt2$, a pentagonal
cases discussed further in Section 5.9.)
  The process terminates
because the number of sides to the polygon drops at every step. A
simple recursive computer procedure runs through all possible ways
the subregion might break into pieces and checks that the
tcc-bound gives Inequalities \refz (5.5.1) and \refz (5.5.2).  The
same argument works if there is a special simplex that overlaps
each of the other special simplices in the subcluster.

When $n=6,7$ and there are two nonoverlapping special simplices,
a similar argument
can be applied.
 Remove both specials by cutting along the diagonals.
  Then deform both diagonals to length
$2\sqrt{2}$, taking into account the possible ways that the subregion
can break into pieces in the process.  In every case the bounds
(5.5.1) and (5.5.2)
are satisfied.

There are a number of situations that arise that escape this generic
argument and were analyzed individually.
These include the cases involving
more than two special simplices over a given subregion,
two special simplices over a pentagon, or a special simplex over a
quadrilateral.
  Also, the deformation lemmas
are insufficient to bring all of the edges between
adjacent corners to one of the three standard lengths $2,2.51,2\sqrt{2}$
for certain triangular and quadrilateral regions.  These are treated
individually.

The next few sections describe the cases treated individually.
  The cases not
mentioned in the sections that follow fall within the generic
procedure just described.

\subhead 5.7 Triangles \endsubhead

With triangular subregions, there is no need to use any of the
deformation arguments because the dimension is already sufficiently
small to apply interval arithmetic directly to obtain our bounds.
There is no need for the tcc-bound approximations.

Flat quarters and simplices of type $S_A$ are treated in $\A_{16}$.
Other simplices are scored by the truncated Voronoi function.
We break the edges between corners into the cases
$[2,2.51)$, $[2.51,2\sqrt{2})$, $[2\sqrt{2},3.2]$. 
Let $k_0$, $k_1$, and $k_2$, with $k_0+k_1+k_2=3$, be the number
of edges  in the respective intervals.

If $k_2=0$, we can improve the penalties,
    $$\pi_\tau = \pi_\sigma=0.$$
To see this, first we observe that there can be no $\Sminus$
or $\Sfour$ configurations.  By placing $\ge3$ quarters
around an
upright diagonal, if the subregion is triangular, the upright
diagonal becomes surrounded by anchored simplices, a case deferred
until Section 5.11.

If $k_0=k_1=k_2=1$, we can take
$\pi'_\tau= \xiG+2\xiV+0.0114=0.034052$.
A few cases are needed to justify this constant.
If there are no $\Sminus$ configurations, $\pi'_\tau$ is
at most
$$
\align
&[\xiG + 2 \xiV +\xikG ]3/4 < 0.0254,\\
\hbox{or\quad }&[\xiG+2\xiV+\xikG]2/4 + 0.008/3 < 0.0254
\endalign
$$
If there are at most two edges
    in the subregion coming from an $\Sminus$ configuration,
    $$(\xiG+2\xiV+0.0114)2/3 + 0.008/3 < 0.0254.$$
If three edges come from an $\Sminus$ configuration, we get
$0.034052$.
To get somewhat sharper bounds, we consider how the edge $k_2$
was formed.  If it is obtained by deformation from an edge
in the standard region of length $\ge3.2$, then it becomes a distinguished
edge when the length drops to $3.2$.  If the edge in the
standard region already has length $\le3.2$, then it is distinguished
before the deformation process begins, so that the subregion
can be treated in isolation from the other subregions.
We conclude that when $\pi'_\tau=0.034052$ we can take
$y_4\ge2.6$ or $y_5=3.2$ (Remark 3.9).

The bounds (5.5.1) and (5.5.2) now follow from $\A_{17}$ and $\A_{18}$.

\subhead 5.8 Quadrilaterals \endsubhead

We introduce some notation for the heights and edge lengths of a
convex polygon.  The heights will generally
be $2$ or $2.51$, the edge lengths
between consecutive corners will generally
be $2$, $2.51$, or $2\sqrt{2}$.  We
represent the edge lengths by a vector
    $$(a_1,b_1,a_2,b_2,\ldots,a_n,b_n),$$
if the corners of an $n$-gon, ordered cyclically have heights
$a_i$ and if the edge length between corner $i$ and $i+1$ is
$b_i$.  We say two vectors are equivalent if they are related by a
different cyclic ordering on the corners of the polygon, that is,
by the action of the dihedral group.

The vector of a polygon with a special simplex is equivalent to
one of the form $(2,2,a_2,2,2,\ldots)$.  If $a_2=2.51$, then
what we have is necessarily special (Section 5.6).
However, if $a_2=2$, it is possible
for the edge opposite $a_2$ to have length greater than $3.2$.

Turning to quadrilateral regions, we use tcc scoring if both
diagonals are greater than 3.2.   Suppose that both diagonals
are between $[2\sqrt{2},3.2]$, creating a pair of overlapping
special simplices.  The deformation lemma requires a diagonal
longer than $3.2$, so although we can bring the quadrilateral
to the form
$$(a_1,2,2,2,2,2,a_4,b_4),$$
the edges $a_1,a_4,b_4$ and the diagonal vary continuously (see
$\A_{13}$).
By $\A_{19}$, we have bounds on the score
$$
\align
\tau_0 &> 0.235, \quad \vor_0 < -0.075,
                \hbox{ if } b_4\in[2.51,2\sqrt{2}],\\
\tau_0 &> 0.3109, \quad \vor_0 < -0.137,
                \hbox{ if } b_4\in[2\sqrt{2},3.2],\\
\endalign
$$
We have $D(4,1)=0.2052$, $Z(4,1)=-0.05705$.
When $b_4\in[2.51,2\sqrt{2}]$, we can take
$\pi_\tau=\pi_\sigma=0$. (We are excluding loops here.)
When $b_4\in[2\sqrt{2},3.2]$, we can take
    $$
    \align
    \pi_\tau &= \maxpi+ 0.0066, \\
    \pi_\sigma &= 0.008 (5/3)+ 0.009. \\
    \endalign
    $$
It follows that the Inequalities \refz (5.5.1) and \refz (5.5.2)
are satisfied.

Suppose that one diagonal has length $[2\sqrt{2},3.2]$ and the
other has length at least $3.2$.  The quadrilateral is represented
by the vector
$$(2,2,a_2,2,2,b_3,a_4,b_4).$$
The hypotheses of the deformation lemma hold,
so that $a_i\in\{2,2.51\}$ and
$b_j\in\{2,2.51,2\sqrt2\}$.
To avoid quad clusters, we assume $b_4\ge\max(b_3,2.51)$.
These are one-dimensional
with a diagonal of length $[2\sqrt{2},3.2]$ as parameter.
The required verifications appear in $\A_{20}$.

\subhead 5.9 Pentagons \endsubhead

Some extra comments are needed when there is a special simplex.
The general argument outlined above removes the special, leaving
a quadrilateral.  The quadrilateral is deformed, bringing the
edge that was the diagonal of the special to $2\sqrt{2}$.
This section discusses how this argument might break down.

Suppose first that there is a special and that both diagonals
on the resulting quadrilateral are at least 3.2.  We can deform
using either diagonal, keeping both diagonals at least 3.2.
The argument breaks down if both diagonals drop to 3.2 before
the edge of the special reaches $2\sqrt{2}$ and both diagonals
of the quadrilateral lie on specials.
When this happens,
 the quadrilateral has the form
$$(2,2,2,2,2,2,2,b_4),$$
where $b_4$ is the edge originally on the special simplex.  If
both diagonals are 3.2, this is rigid, with $b_4= 3.12$.
  We find its score to be
$$
\align
&\vor_0(S(2,2,2,b_4,3.2,2))+\vor_0(S(2,2,2,3.2,2,2))+0.0461<-0.205,\\
&\tau_0(S(2,2,2,b_4,3.2,2))+\tau_0(S(2,2,2,3.2,2,2))2> 0.4645.\\
\endalign
$$
So the Inequalities \refz (5.5.1) and \refz (5.5.2) hold easily.

If there is a special and there is a diagonal on
resulting quadrilateral $\le3.2$, we have two nonoverlapping
specials.  It has the form
$$(2,2,a_2,2,2,2,a_4,2,2,b_5).$$
The edges $a_2$ and $a_4$ lie on the special.  If $b_5>2$,
cut away one of the special simplices.  What is left can be
reduced
to a triangle, or a quadrilateral case treated in $\A_{20}$.
Assume $b_5=2$.  We have a pentagonal standard region.
We may assume that there is no $\Sfour$ or $\Sminus$ configuration,
for otherwise Theorem 4.4 follows trivially from
the bounds in Section 2.
  A pentagon can then have at most
$\Splus$ for a penalty of $0.008$.

If $a_2=2.51$ or $a_4=2.51$, we
again remove a special simplex and produce triangles, quadrilaterals,
or the special cases in $\A_{20}$.
We may impose the condition $a_2=a_4=b_5=2$.
We score this full pentagonal arrangement in $\A_{21}$, using
the edge lengths of the two diagonals of the specials as variables.
The inequalities follow.

\subhead 5.10 Hexagons and heptagons \endsubhead

We turn to hexagons.
There may be three specials whose diagonals do not cross.  Such
a subcluster is represented by the vector
$$(2,2,a_2,2,2,2,a_4,2,2,2,a_6,2).$$
The heights $a_{2i}$ are $2$ or $2.51$.  Draw the diagonals between
corners $1$, $3$, and $5$.  This is a three-dimensional configuration,
determined by the lengths of the three diagonals.  The required
bound follows from $\A_{21}$.

There is one case with a special simplex that
did not satisfy the generic computer-checked inequalities for
what is to be squandered.  Its vector is
    $$(a_1,2,2,2,2,2,2,b_4,2,2,2,2),$$
with $a_1=b_4=2.51$.
A vertex of the special simplex has height $a_1=2.51$ and
all other corners have height $2$.  The subregion is a hexagon
with one edge longer than $2$.  We have $D(6,1)= 0.48414$.
This is certainly obtained if the subregion contains the configuration
$\Sminus$ squandering $0.5606$.
But if this configuration does
not appear, we can decrease $\pi_\tau$ to
    $0.03344 + (2/3) 0.008$, a constant coming from $\Sfour$
 in Section 4.7.  With this smaller penalty the
inequality is satisfied.

Now turn to heptagons.
The bound $2\pi$  on the perimeter of the polygon, eliminates
all but one equivalence class of vectors associated with a polygon
that has two or more potentially specials simplices. The vector is
    $$(2,2,a_2,2,2,2,a_4,2,2,2,a_6,2,a_7,2),$$
$a_2=a_4=a_6=a_7=2.51$.
In other words, the edges between adjacent corners are $2$ and
four heights are $2.51$. There are two specials.
  This case is treated by the procedure outlined for
subregions with two specials whose diagonals do not cross.

\subhead 5.11 Loops \endsubhead

We now return to a collection of anchored simplices that surround
the upright diagonal.  This is the last case needed to complete
the proof of Theorem 4.5.
There are four or five anchored simplices around the upright diagonal.
$\A_2$--$\A_7$ give a list of linear inequalities satisfied by the
anchored simplices, broken up according to type: upright, type $S_C$,
opposite edge $>3.2$, etc.   The anchored simplices are related
by the constraint that the sum of the dihedral angles around
the upright diagonal is $2\pi$.  We run a linear program in each
case based on these linear inequalities, subject to this constraint
to obtain bounds on the score and what is squandered by the
anchored simplices.

When the edge opposite the diagonal of an anchored simplex
has length $\in[2\sqrt{2},3.2]$
and the simplex adjacent to the anchored simplex across that
edge is a special simplex, we use the inequalities $\A_{22}$ and
$\A_{23}$ that run parallel to $\A_4$ and $\A_5$.  It is not necessary
to run separate linear programs for these.  It is enough to observe
that the constants for what is squandered improve on those
from $\A_4$ by at least $0.06445$ and that the constants for
the score in $\A_{22}$ differ with those of $\A_4$ by no more
than $0.009$.

When the dihedral angle of an anchored simplex is greater than
$2.46$, the simplex is dropped, and the remaining anchored simplices
are subject to the constraint that their dihedral angles sum
to at most $2\pi-2.46$.  There can not be an anchored simplex
with dihedral angle greater than $2.46$ when there are five
anchors: $2.46+4 (0.956)>2\pi$.
There cannot be two anchored simplices
with dihedral angle greater than $2.46$: $2(2.46+0.956)>2\pi$
($\A_8$).

The following table summarizes the linear
programming results.

$$
\matrix
(n,k)   &   \DLP(n,k) & D(n,k)      &\ZLP(n,k)  &Z(n,k)\\
(4,0)   &   0.1362  &   0.1317  &   0   &   0\\
(4,1)   &   0.208   &   0.20528 &-0.0536&   -0.05709\\
(4,2)   &   0.3992  &   0.27886 &-0.2   &   -0.11418\\
(4,3)   &  0.6467   &   0.35244 &-0.424 &   -0.17127\\
(5,0)   &   0.3665  &   0.27113 &-0.157 &   -0.05704\\
(5,1)   &  0.5941   &   0.34471 &-0.376 &   -0.11413\\
(5,\ge2)&  0.9706   &  \squander    &*          &   *
\endmatrix
$$

The bound for $D(4,0)$ comes from III.4.1.11.
A few more comments are needed for $Z(4,1)$.  Let $S=S(y_1,\ldots,y_6)$
be the anchored simplex that is not a quarter.  If $y_4\ge2\sqrt2$
or $\dih(S)\ge 2.2$, the linear programming bound is $<Z(4,1)$.
With this, if $y_1\le 2.75$, we have $\sigma(S) < Z(4,1)$
by $\A_{12}$.  But if $y_1\ge2.75$, the 3 upright quarters along the
upright diagonal satisfy
    $$\nu< -0.3429+0.24573\dih.$$
With this stronger inequality, the linear programming bound
becomes $<Z(4,1)$.
This completes the proof of Theorem 4.4.
\qed

\bigskip
\subhead 5.12 Some final estimates \endsubhead
\medskip

Recall that Section \refz 4.4 defines an integer $n(R)$ that is
equal to the number of sides if the region is a polygon.  Recall
that if the dihedral angle along an edge of a standard cluster is
at most $1.32$, then there is a flat quarter along that edge
(Lemma \refz 3.11.4).

\proclaim{Lemma 5.12.1}  Let $R$ be an exceptional cluster with a
dihedral angle $\le1.32$ at a vertex $v$. Then $R$ squanders
$>t_n+1.47\,\pt$, where $n=n(R)$.
\endproclaim

\demo{Proof} In most cases we establish the stronger bound
$t_n+1.5\,\pt$. In the proof of Theorem \refz 4.4, we erase all
upright diagonals, except those completely surrounded by anchored
simplices.  The contribution to $t_n$ from the flat quarter $Q$ at
$v$ in that proof is $D(3,1)$ (Sections \refz 4.5 and \refz
5.5.1). Note that
 $\epsilon_\tau(Q)=0$ here because there
are no deformations. If we replace $D(3,1)$ with $3.07\,\pt$ from
Lemma \refz 3.11.4, then we obtain the bound. Now suppose the
upright diagonal is completely surrounded by anchored simplices.
  Analyzing the constants of Section \refz 5.11,
we see that $\DLP(n,k)-D(n,k)>1.5\,\pt$. except when
$(n,k)=(4,1)$.

Here we have four anchored simplices around an upright diagonal.
Three of them are quarters.  We erase and take a penalty. Two
possibilities arise.  If the upright diagonal is enclosed over the
flat quarter, its height is $\ge2.6$ by geometric considerations
and the top face of the flat quarter has circumradius at least
$\sqrt2$.  The penalty is $2\xiG' + \xiV$, so the bound holds by
the last statement of Lemma \refz 3.11.4.

If, on the other hand, the upright diagonal is not enclosed over
the flat diagonal, the penalty is $ \xiG + 2\xiV$.  In this case,
we obtain the weaker bound $1.47\,\pt+t_n$:
    $$3.07\,\pt > D(3,1) + 1.47\,\pt +\xiG+2\xiV.$$
\qed\enddemo

\proclaim{Remark}  If there are $r$ nonadjacent vertices with
dihedral angles $\le1.32$, we find that $R$ squanders
$t_n+r(1.47)\,\pt$.
\endproclaim

In fact, in the proof of the lemma, each $D(3,1)$ is replaced with
$3.07\,\pt$ from Lemma \refz 3.11.4.  The only questionable case
occurs when two or more of the vertices are anchors of the same
upright diagonal (a loop). Referring to Section \refz 5.11, we
have the following observations about various contexts.

{ \noindent
\parskip=0pt
\parindent=0pt
\hbox{}

$(4,1)$ can mask only one flat quarter and it is treated in the
lemma.

$(4,2)$ can mask only one flat quarter and
$\DLP(4,2)-D(4,2)>1.47\,\pt$.

$(4,3)$ cannot mask any flat quarters.

$(5,0)$ can mask two flat quarters.  Erase the five upright
quarters, and
    take a penalty $4\xiV+\xiG$.  We get
    $$D(3,2)+2(3.07)\,\pt > t_5+4\xiV+\xiG+2(1.47)\,\pt.$$

$(5,1)$ can mask two flat quarters, and
$\DLP(5,1)-D(5,1)>2(1.47)\,\pt$.

}

\proclaim{Lemma 5.12.2} Any pentagon with a dihedral angle less
than $1.32$ squanders at least $5.66\,\pt$.
\endproclaim

\demo{Proof} To obtain the bound $5.66\,\pt$, we argue as follows.
If there are five anchored simplices surrounding a vertex, we have
the bound by Table \refz 5.11.  If the configuration $\Sfour$ or
$\Sminus$ occurs, we squander at least
    $0.4 > 5.66\,\pt$
(Sections \refz 3.8 and \refz 3.7). So if there are any upright
diagonals in the pentagon that carry a penalty, we may assume they
have four anchors.  If there are no penalties, Lemma \refz 3.11.4
gives
    $3.07\,pt+D(4,1)>5.66\,\pt$.
We do not need to deal with penalties from $\Splus$ in the score
of the flat quarter at $v$ because all penalties from a flat
quarter are applied to the adjacent subregion (see Section \refz
5.5 and Lemma \refz 3.9.2). The only remaining possibility is four
anchored simplices surrounding an upright diagonal.  Unless there
are three upright quarters, the bound follows from Section \refz
5.11.  If there are three upright quarters, erasing gives penalty
$3\xiG$, and
    $3.07\,\pt+D(4,1)-3\xiG>5.66\,\pt$.
This proves the lemma for two pentagons and a quadrilateral.
\qed\enddemo


\vfill\eject
\parskip=0.2\parskip
\head References\endhead

\noindent
[F] S. Ferguson, T. Hales, A Formulation of the Kepler
    Conjecture, preprint

\noindent
[HM] T. Hales, S. McLaughlin, A  Proof of the Dodecahedral Conjecture,
    preprint

\noindent
[I] Thomas C. Hales, Sphere Packings I,
    Discrete and Computational Geometry, 17 (1997), 1-51.

\noindent
[II] Thomas C. Hales, Sphere Packings II,
    Discrete and Computational Geometry, 18 (1997), 135-149.

\noindent
[III] Thomas C. Hales, Sphere Packings III, preprint.

\noindent
[V] S. Ferguson, Sphere Packings V, thesis, University of Michigan,
    1997.

\noindent
[H1] Thomas C. Hales, Packings, \hfill\break
    {\rm http://www.math.lsa.umich.edu/\~%
    \relax hales/packings.html}

\noindent
[H2] Thomas C. Hales, Remarks on the Density of Sphere Packings,
        Combinatorica, 13 (2) (1993) 181-197.

\heads{Appendix 1}
\vfill\eject
\head Appendix 1. Inequalities\endhead

\footnote""{-- {\it Sphere Packings IV} --
printed \rm\today }

Let $\octavor_0(Q)=0.5(\vor_0(Q)+\vor_0(\hat Q))$, and
$\octavor(Q) = 0.5(\vor(Q)+\vor(\hat Q))$.
We let $\tau_\nu$, $\tau_V$, $\tau_{0}$, and $\tau_\Gamma$ be
$-f+\sol\zeta\pt$, where $f=\nu$, $\vor$, $\vor_0$, and $\Gamma$,
respectively.

Each inequality is accompanied by one or more reference numbers.
These identification numbers are needed to find further details
about these calculations in \cite{H1}.
These inequalities were checked numerically
before they were rigorously established,
using a nonlinear
optimization package.
I thank the University of Maryland
for this software.
\footnote"*"{\tt \quad www.isr.umd.edu/Labs/CACSE/FSQP/fsqp.html}

  Edge lengths whose bounds are not specified are assumed
to be between 2 and $2.51$.

Most of the interval calculations in this appendix were completed
by Samuel Ferguson.  His calculations are marked with a dagger (\dag).

\def\refno#1{\hbox{}\nobreak\hfill {\tt (#1)}}

\parindent=0pt

\subhead Section $\A_1$\endsubhead
    $\beta_\psi$ is defined in Section \refz 2.8.

 1: $\beta_\psi(y_1,y_3,y_5) <\dih_3(S)$,
if $y_2,y_3\in[2,2.23]$, $y_4\in[2.77,2\sqrt{2}]$, $\cos\psi=y_1/2.77$.
(We may assume $y_6=2$.)
    \refno{757995764\dag}

 2: $\beta_\psi(y_1,y_3,y_5) <\dih_3(S)$,
provided  $y_4=3.2$, $y_5=2.51$, $y_6=2$, $\cos\psi=y_1/2.51$.
\refno{735258244\dag}

 3: $\beta_\psi(y_1,y_2,y_6)<\dih_2(S)$,
    if $y_4\in[2,3.2]$, $y_5=y_6=2.51$, $\psi=\arc(y_1,t_0,1.6)$.
\refno{343330051\dag}

4:  $\beta_\psi(y_1,y_2,y_6)<\dih_2(S)$,
    if $y_4\in[2,3.2]$, $y_5=3.2$, $y_6=2$, $y_1\in[2.2,2.51]$,
    $\psi=\arc(y_1,t_0,1.6)$.
\refno{49446087\dag}

5: $\dih(R(y_2/2,\eta_{126},y_1/(2\cos\psi)))<\dih_2(S)$, if
$y_1\in[2,2.2]$, $y_3=2.51$, $y_4=3.2$, $y_5=3.2$, $y_6=2$,
$\psi=\arc(y_1,t_0,1.6)$. \refno{799187442\dag}

6:  $\vor(Q,1.385)<0.00005$, if $y_4\in[2.77,2\sqrt{2}]$,
        and $\eta_{456}\ge\sqrt2$.
    \refno{275706375}

7:  $\vor(Q,1.385)<0.00005$, if $y_4\in[2.77,2\sqrt{2}]$,
        and $\eta_{234}\ge\sqrt2$.
    \refno{324536936}

8:  $\tau_V(Q,1.385)>0.0682$, if $y_4\in[2.77,2\sqrt{2}]$,
        and $\eta_{456}\ge\sqrt2$.
    \refno{983547118}

9:  $\tau_V(Q,1.385)>0.0682$, if $y_4\in[2.77,2\sqrt{2}]$,
        and $\eta_{234}\ge\sqrt2$.
    \refno{206278009}

\subhead Section $\A_2$\endsubhead
In Inequalities $\A_2$ and $\A_3$,
the domain is the set of upright quarters.
The dihedral angle is measured along the diagonal.

1: $\nu < -4.3223 + 4.10113\dih$.
    \refno{413688580}

2: $\nu < -0.9871 + 0.80449\dih$
    \refno{805296510}

3: $\nu < -0.8756 + 0.70186\dih$
    \refno{136610219}

4: $\nu < -0.3404 + 0.24573\dih$
    \refno{379204810}

5: $\nu < -0.0024 + 0.00154\dih$
    \refno{878731435}

6: $\nu < 0.1196 - 0.07611\dih$
    \refno{891740103}

\subhead Section $\A_3$\endsubhead

1: $-\tau_\nu < -4.42873 + 4.16523\dih$
    \refno{334002329}

2: $-\tau_\nu < -1.01104 + 0.78701\dih$
    \refno{883139937}

3: $-\tau_\nu < -0.99937 + 0.77627\dih$
    \refno{507989176}

4: $-\tau_\nu < -0.34877 + 0.21916\dih$
    \refno{244435805}

5: $-\tau_\nu < -0.11434 + 0.05107\dih$
    \refno{930176500}

6: $-\tau_\nu < 0.07749 - 0.07106\dih$
    \refno{815681339}

\subhead Section $\A_4$\endsubhead

In $\A_4$ and $\A_5$,  $y_1\in[2.51,2\sqrt{2}]$ and
$y_4\in[2.51,2\sqrt{2}]$ and $\dih<2.46$.
Let $\vor_x = \vor$ if the simplex is of type $C$ and
$\vor_x=\vor_0$ otherwise.

1:  $\vor_x < -3.421 + 2.28501 \dih$
    \refno{649592321}

2: $\vor_x < -2.616 + 1.67382 \dih$
    \refno{600996944}

3: $\vor_x < -1.4486 + 0.8285 \dih$
    \refno{70667639}

4: $\vor_x < -0.79 + 0.390925 \dih$
    \refno{99182343}

5: $\vor_x < -0.3088 + 0.12012 \dih$
    \refno{578762805}

6: $\vor_x < -0.1558 + 0.0501 \dih$
    \refno{557125557}

\subhead Section $\A_5$\endsubhead
Set $\tau_x =\sol\zeta\pt - \vor_x$.

1: $-\tau_x < -3.3407 + 2.1747 \dih$
    \refno{719735900}

2: $-\tau_x < -2.945 + 1.87427 \dih$
    \refno{359616783}

3: $-\tau_x < -1.5035 + 0.83046 \dih$
    \refno{440833181}

4: $-\tau_x < -1.0009 + 0.48263 \dih$
    \refno{578578364}

5: $-\tau_x < -0.7787 + 0.34833 \dih$
    \refno{327398152}

6: $-\tau_x < -0.4475 + 0.1694 \dih$
    \refno{314861952}

7: $-\tau_x < -0.2568 + 0.0822 \dih$
    \refno{234753056}

\subhead Section $\A_6$\endsubhead
In the Inequalities $\A_6$ and $\A_7$, we assume $y_1\in[2.51,2\sqrt{2}]$,
$y_4\in[2\sqrt{2},3.2]$, and $\dih<2.46$.

1:  $\vor_0 < -3.58 + 2.28501 \dih$
    \refno{555481748}

2: $\vor_0 < -2.715 + 1.67382 \dih$
    \refno{615152889}

3: $\vor_0 < -1.517 + 0.8285 \dih$
    \refno{647971645}

4: $\vor_0 < -0.858 + 0.390925 \dih$
    \refno{516606403}

5: $\vor_0 < -0.358 + 0.12012 \dih$
    \refno{690552204}

6: $\vor_0 < -0.186 + 0.0501 \dih$
    \refno{852763473}

\subhead Section $\A_7$\endsubhead The assumptions are as in
$\A_6$.

1: $-\tau_0 < -3.48 + 2.1747 \dih$
    \refno{679673664}

2: $-\tau_0 < -3.06 + 1.87427 \dih$
    \refno{926514235}

3: $-\tau_0 < -1.58 + 0.83046 \dih$
    \refno{459744700}

4: $-\tau_0 < -1.06 + 0.48263 \dih$
    \refno{79400832}

5: $-\tau_0 < -0.83 + 0.34833 \dih$
    \refno{277388353}

6: $-\tau_0 < -0.50 + 0.1694 \dih$
    \refno{839852751}

7: $-\tau_0 < -0.29 + 0.0822 \dih$
    \refno{787458652}

\subhead Section $\A_8$\endsubhead In all these except {\tt
(125103581)} and {\tt (504968542)}, the signs of all the partials
except in the $x_1$ variable are easily determined by the methods
of Section I.8.  In this way, they become optimizations in one
variable.

1:    $\dih>1.23$ if $y_1\in [2.51,2\sqrt{2}]$, and $y_4\ge 2.51$.
    \refno{499014780}

2:    $\dih > 1.4167$, if $y_1\in[2.51,2\sqrt{2}]$, and
    $y_4\ge2\sqrt{2}$.
    \refno{901845849}

3: $\dih> 1.65$ if $y_1\in[2.51,2\sqrt{2}]$, $y_4\ge3.2$
    \refno{410091263}

4: $\dih> 0.956$ if $y_1\in[2.51,2\sqrt{2}]$, $y_4\ge 2$.
    \refno{125103581}

5: $\dih > 0.28$, if $y_1\in[2.51,2\sqrt{2}]$, $y_4\ge2$,
    $y_5\in[2,2\sqrt{2}]$.
    \refno{504968542}

6: $\dih> 1.714$, if $y_1\in[2.7,2\sqrt{2}]$, $y_4\ge3.2$
    \refno{770716154}

7: $\dih> 1.714$, if $y_1\in[2.51,2.7]$, $y_4\ge3.2$,
    $y_2\in[2,2.25]$
    \refno{666090270}

8:  $\dih < 2.184$, if  $y_1\in[2.51,2\sqrt{2}]$.
    (This one was simple enough to do without interval arithmetic.)
    \refno{971555266}

\subhead Section $\A_9$\dag\endsubhead
$\kappa(S)$ is defined in Section \refz 3.3.

1: $\kappa< -0.003521$,
    $y_1\in[2.696,2\sqrt{2}]$, $y_2,y_6\in[2.45,2.51]$,
    $y_4\ge 2.77$,
    \refno{956875054}

2: $\kappa < -0.017$,
    if $y_1\in[2.51,2.696]$, $y_4\in[2.77,2\sqrt{2}]$, $\eta_{234}\ge \sqrt{2}$.
    \refno{664200787}

3: $\kappa < -0.017$,
    if $y_1\in[2.51,2.696]$, $y_4\in[2.77,2\sqrt{2}]$,
        $\eta_{456}\ge \sqrt{2}$.
    \refno{390273147}

4: $\kappa< -0.02274 = \xikG-\xiG'$, if
    $y_1\in[2.57,2\sqrt{2}]$, $y_4\ge 3.2$,
    $\Delta\ge0$.
By monotonicity we may assume $y_4=3.2$.
    \refno{654422246}

5: $\kappa< \xik = -0.029$, if
    $y_1\in[2.51,2.57]$, $y_4\ge 3.2$,
    $\Delta\ge0$.
By monotonicity we may assume $y_4=3.2$.
    \refno{366536370}

6:  $\kappa< -0.03883$,
    if $y_1\in[2.51,2.57]$,
    $y_2,y_3,y_5,y_6\in[2,2.25]$,
    $y_4\ge3.2$, $\Delta\ge0$.
By monotonicity we may assume $y_4=3.2$.
\refno{62532125}

7:  $\kappa< -0.0325$,
    if $y_1\in[2.51,2.57]$,
    $y_2,y_3,y_5\in[2,2.25]$,
    $y_4\ge3.2$,
    $\Delta\ge0$.
By monotonicity we may assume $y_4=3.2$.
    \refno{370631902}

\subhead Section $\A_{10}$\endsubhead

1:  $\Gamma< \octavor_0$, if
    $y_1\in[2.696,2\sqrt{2}]$.
    \refno{214637273}

2: $\Gamma < \octavor_0 + 0.01561$,
    if $y_1\in[2.51,2\sqrt{2}]$.
    \refno{751772680}

3: $\Gamma < \octavor_0 + 0.00935$, if
    $y_1\in[2.57,2\sqrt{2}]$.
    \refno{366146051}

4: $\Gamma<\octavor_0+0.00928$,
    if $y_1\in[2.51,2.57]$,
    $y_2\in [2.25,2.51]$.
    \refno{675766140}

5: $\Gamma<\octavor_0$,
    if $y_1\in[2.51,2.57]$,
    $y_2,y_6\in [2.25,2.51]$.
    \refno{520734758}

\subhead Section $\A_{11}$\endsubhead

1: $\octavor<\octavor_0$, if $y_1\in[2.696,2\sqrt{2}]$,
$y_2,y_3\in[2,2.45]$.
    \refno{378432183}

2: $\octavor< \octavor_0$, if $y_1\in[2.696,2\sqrt{2}]$,
    $y_2,y_5\in[2.45,2.51]$.
    \refno{572206659}

3:  $\vor<\vor_0 + 0.003521$, if $y_1\in[2.51,2\sqrt{2}]$.
    \refno{310679005}

4: $\vor <\vor_0 - 0.003521$, if
    $y_1\in[2.696,2\sqrt{2}]$, $y_2,y_6\in[2.45,2.51]$,
    $y_4\in[2.51,2.77]$.
    \refno{284970880}

5: $\vor < \vor_0 -0.009$.
    if $y_1\in[2.51,2.696]$, $y_4\in[2.51,2\sqrt2]$.
    \refno{972111620}

6: $\octavor < \octavor_0$,
    if $y_1\in[2.51,2.57]$,
    $\eta_{126}\ge\sqrt2$.
    \refno{875762896}

7: $\octavor<\octavor_0 - 0.004131$,
    if $y_1\in[2.51,2\sqrt{2}]$,
    $\eta_{126}\le\sqrt2$,
    $\eta_{135}\ge\sqrt2$,
    $y_3\le2.2$.
    \refno{385332676}

\subhead Section $\A_{12}$\endsubhead

1: $\tau_V(S) > 0.13 + 0.2 (\dih(S)-\pi/2)$,
    if $y_1,y_2\in[2.51,2\sqrt{2}]$, and $\eta_{126}(S)\le\sqrt{2}$.
    \refno{970291025\dag}

2: $\tau_V(S,\sqrt{2}) > 0.13 + 0.2 (\dih(S)-\pi/2)$,
    if $y_1,y_2\in[2.51,2\sqrt{2}]$, and $\eta_{126}(S)\ge\sqrt{2}$.
    \refno{524345535\dag}

3: $\nu < -0.3429 + 0.24573\dih$, for upright quarters with
    $y_1\in[2.75,2\sqrt{2}]$.
    \refno{812894433}

4: $\vor_x < -0.0571$, for anchored simplices with
    $y_4\in[2.51,2\sqrt2]$, $y_1\in[2.51,2.75]$, $\dih<2.2$.
    \refno{404793781}

\subhead Section $\A_{13}$\endsubhead
Inequalities {\tt (74657942)} and {\tt (675901554)}
hold by inspection.  The others are verified in the usual
manner.

1: $\tau_\nu(S)>0.033$, if $S$ is an upright quarter.
    \refno{705592875}

2: $\tau_0(S) > 0.06585- 0.0066$,
    if $S$ is a flat quarter, and $y_4=2\sqrt{2}$.
    \refno{747727191}

3: $\vor_0(S) < 0.009$, if $S$ is a flat quarter, and
    $y_4=2\sqrt{2}$.
    \refno{474496219}

4: $\vor_0(S(2,y_2,y_3,y_4,2,2))<0.0461$, if
    $y_4\in[2\sqrt{2},3.2]$.
    \refno{649551700}

5: $\vor_0(S(2.51,2,y_3,y_4,2,2))\le0$, if
    $y_4\in[2\sqrt{2},3.2]$.
    \refno{74657942}

6: $\vor_0(S(y_1,y_2,2.51,y_4,2,2))<0$, if
    $y_4\in[2\sqrt{2},3.2]$.
    \refno{897129160}

7: $\tau_0(S(2,y_2,y_3,y_4,2,2))>0.014$,  if
    $y_4\in[2\sqrt{2},3.2]$.
    \refno{760840103}

8: $\tau_0(S(2.51,2,2,y_4,2,2))\ge0$,  if $y_4\in[2\sqrt{2},3.2]$.
    \refno{675901554}

9: $\tau_0(S(y_1,y_2,2.51,y_4,2,2))>0.06585$, if
    $y_4\in[2\sqrt{2},3.2]$.
    \refno{712696695}


10: $\nu < \vor_0 + 0.01 (\pi/2-\dih)$, if
    $y_1\in[2.696,2\sqrt2]$. \refno{269048407}

11: $\nu < \vor_0$, if $y_1\in[2.6,2.696]$, $y_4\in[2.1,2.51]$.
    \refno{553285469}

12: $\mu < \vor_0 + 0.0268$, if $y_4\in[2.51,2\sqrt2]$.
    \refno{293389410}

13: $\mu < \vor_0 + 0.02$, if $y_1\in[2,2.17]$,
    $y_4\in[2.51,2\sqrt2]$. \refno{695069283}

14: $\dih>1.32$, if $y_4=2\sqrt{2}$.\refno{814398901}

15: $\hat\tau>3.07\,\pt$, for all flat quarters satisfying
    $\dih\le1.32$.
    \refno{352079526}

16: $\tau_0>3.07\,\pt+\xiV+2\xiG'$, if $y_4\in[2.51,2\sqrt2]$,
    $\eta_{456}\ge\sqrt2$, $\dih\le1.32$.
    \refno{179025673}


\subhead Section $\A_{14}$\dag\endsubhead

$V_i$ is defined in Section \refz 4.9. The function $f$ is defined
in Section \refz 4.13.

1: $V_0 < 0$, if $\Delta\ge0$, $y_4\in[2,y_2+y_3]$,
    $y_5\in[2,3.2]$,
        $y_6\in[y_5,3.2]$.
    \refno{424011442}

2: $V_1 < 0$, if $\Delta\ge0$, $y_4\in[2,y_2+y_3]$,
    $y_5\in[2,3.2]$,
        $y_6\in[y_5,3.2]$.
    \refno{140881233}

3: $V_j + 0.82\sqrt{421}<0$, if $y_5\in[2,2.189]$,
    $y_4\in[2\sqrt{2},3.2]$,
        $y_6\in[2,2.51]$, $\Delta\ge0$, $j=0,1$.
    \refno{601456709}

4: $V_j + 0.82\sqrt{421}<0$, if $y_5\in[2,2.189]$,
    $y_4\in[3.2,y_2+y_3]$,
    $y_6\in[2,3.2]$, $\Delta\ge0$, $j=0,1$.
    \refno{292977281}

5: $V_j + 0.5\sqrt{421}<0$, if $y_5\in[2.189,2.51]$,
$y_4\in[2\sqrt{2},3.2]$,
        $y_5,y_6\in[2,2.51]$, $\Delta\ge0$, $j=0,1$.
    \refno{927286061}

6: $V_j + 0.5\sqrt{421}<0$, if $y_5\in[2.189,3.2]$,
$y_4\in[3.2,y_2+y_3]$,
        $y_5,y_6\in[2,3.2]$, $\Delta\ge0$, $j=0,1$.
    \refno{340409511}

7: $\Delta<421$, if $y_4\in[2\sqrt{2},y_2+y_3]$,
$y_5,y_6\in[2,3.2]$,
    $\eta(x_1,x_3, x_5)\le t_0$.
    \refno{727498658}

8: $-4\doct u_{135} \partial/\partial x_5
(\quo(R_{135})+\quo(R_{315}))< 0.82$.
    \refno{484314425}

9: $-4\doct u_{135} \partial/\partial x_5
(\quo(R_{135})+\quo(R_{315}))< 0.5$,
    if $y_5\in[2.189,2.51]$.
    \refno{440223030}

10: $f(y_1,y_2)\ge 0.887$, $\lambda=1.945$, $y_1,y_2\in[2,2.51]$.
    \refno{115756648}

\subhead Section $\A_{15}$\dag\endsubhead
Let $D^if_j = \partial^i f_j(S)/\partial x_1^i$,
$f_0=\vor_0$, $f_1=-\tau_0$, as in Section 5.1.

1: $D^2f_i>0$ if $Df_i=0$, if $\Delta\ge0$, $y_4\ge2\sqrt2$,
    $y_5=2$, $y_6=2$, $y_4\le y_2+y_3,y_5+y_6$,
    $i=0,1$.
    \refno{329882546}

2: $D^2f_i>0$ if $Df_i=0$, if $\Delta\ge0$, $y_4\ge2\sqrt2$,
$y_5=2$, $y_6=2.51$,
    $y_4\le y_2+y_3,y_5+y_6$,
    $i=0,1$.
    \refno{427688691}

3: $D^2f_i>0$ if $Df_i=0$,
 if $\Delta\ge0$, $y_4\ge2\sqrt2$, $y_5=2$, $y_6=2\sqrt{2}$,
    $y_4\le y_2+y_3,y_5+y_6$,
    $i=0,1$.
    \refno{562103670}

4: $D^2f_i>0$ if $Df_i=0$,
 if $\Delta\ge0$, $y_4\ge2\sqrt2$, $y_5=2.51$, $y_6=2.51$,
    $y_4\le y_2+y_3,y_5+y_6$,
    $i=0,1$.
    \refno{564506426}

5: $D^2f_i>0$ if $Df_i=0$,
 if $\Delta\ge0$, $y_4\ge2\sqrt2$, $y_5=2.51$, $y_6=2\sqrt{2}$,
    $y_4\le y_2+y_3,y_5+y_6$,
    $i=0,1$.
    \refno{288224597}

6: $D^2f_i>0$ if $Df_i=0$,
    $\Delta\ge0$, $y_4\ge2\sqrt2$, $y_5=2\sqrt{2}$, $y_6=2\sqrt{2}$,
    $y_4\le y_2+y_3,y_5+y_6$.
    $i=0,1$.
    \refno{979916330,749968927}

\subhead Section $\A_{16}$\endsubhead
Recall $D(3,2)=0.13943$, $Z(3,2)=-0.05714$, $D(3,1)=0.06585$.
Some of these follow from known results.  See II.4.5.1, F.3.13.1,
F.3.13.3, F.3.13.4.  The case $\vor\le0$ of the inequality
$\sigma\le0$ for flat quarters follows by Rogers's monotonicity
lemma I.8.6.2 and F.3.13.1, because the circumradius of the
flat quarter is at least $\sqrt2$ when the analytic Voronoi function
is used.  We also use that $\vor(R(1,\eta(2,2,2)\sqrt2))=0$.

1: $\tilde\tau(S) > 0.06585$, if $S$ is a flat quarter and
$\tilde\tau(S)$ is
    any of the functions for flat quarters in Section 3.10,
    other than $\tau(S,1.385)$, which is treated in $\A_1$.
    \refno{695180203}

2: $\tilde\sigma(S) \le 0$, if $S$ is a flat quarter and
$\tilde\sigma(S)$ is
    any of the functions for flat quarters in Section 3.10,
    other than $\vor(S,1.385)$, which is treated in $\A_1$.
    \refno{690626704}

3: $\vor(S) < Z(3,2)$, for simplices $S$ of type $S_A$.
    \refno{807023313}

4: $\tau_V(S)> 0.13943$, for simplices $S$ of type $S_A$.
    \refno{590577214}

5: $\vor_0(S) < Z(3,2)$, if $y_4,y_5\in[2.51,2\sqrt2]$, and the
    simplex $S$ is not of type $S_A$.
    \refno{949210508}

6: $\tau_0(S) > 0.13943$, if $y_4,y_5\in[2.51,2\sqrt2]$, and the
    simplex $S$ is not of type $S_A$.
    \refno{671961774}

\subhead Section $\A_{17}$\dag\endsubhead
Let $y_4,y_5,y_6\in[2,3.2]$.  Let $k_0$, $k_1$, $k_2$ be the
number of variables in $[2,2.51]$, $[2.51,2\sqrt{2}]$, $[2\sqrt{2},3.2]$,
respectively.(Make the intervals disjoint so that $k_0+k_1+k_2=3$.)
Assume $k_1+2k_2>2$. ($k_1+2k_2=2$ gives special simplices
or cases treated in $\A_{16}$.) We have $3\xiG = 0.04683$.

Set $$\pi'_\tau = \cases 0, & k_2=0,\\
                    0.0254, & k_0=k_1=k_2=1,\\
                    0.04683+(k_0+2k_2-3)0.008/3+k_2(0.0066),& \text{otherwise}.
            \endcases
$$

1: $\tau_0(S) -\pi'_\tau > D(3,k_1+k_2)$, for parameters
$(k_0,k_1,k_2)$ satisfying $k_0+k_1+k_2=3$, $k_1+2k_2>2$.
    \refno{645264496}

2: $\tau_0(S)-0.034052 >D(3,2)$, if $y_4\in[2.6,2\sqrt2]$,
    $y_5\in[2\sqrt2,3.2]$.
    \refno{910154674}

3: $\tau_0(S)-(0.034052+0.0066) >D(3,2)$,
    if $y_6=2,y_4=2.51,y_5=3.2$.
    \refno{877743345}

\subhead Section $\A_{18}$\dag\endsubhead
In the same context as $\A_{17}$,
set

$$\pi'_\sigma = \cases 0, & k_2=0,\\
                    0.009, & k_0=0,k_2=1,\\
                    (k_0+2k_2)0.008/3 + 0.009k_2,& \text{otherwise}.
            \endcases
$$

1: $\vor_0(S) + \pi'_\sigma < Z(3,k_1+k_2)$, for parameters
    $(k_0,k_1,k_2)$ as above.
    \refno{612259047}

\subhead Section $\A_{19}$\dag\endsubhead
Let $Q$ be a quadrilateral subcluster whose edges are described
by the vector
$$(a_1,2,2,2,2,2,a_4,b_4).$$  Assume both diagonals
have lengths in $[2\sqrt2,3.2]$.
$$
\align
\tau_0(Q) &> 0.235 \quad\text{and } \vor_0(Q) < -0.075,
                \hbox{ if } b_4\in[2.51,2\sqrt{2}],\\
\tau_0(Q) &> 0.3109 \quad\text{and } \vor_0(Q) < -0.137,
                \hbox{ if } b_4\in[2\sqrt{2},3.2],\\
\endalign
$$
    \refno{357477295}

\subhead Section $\A_{20}$\dag\endsubhead
Let $Q$ be a quadrilateral subcluster whose edges are described
by the vector $$(2,2,a_2,2,2,b_3,a_4,b_4).$$  Assume
$b_4\ge b_3$, $b_4\in\{2.51,2\sqrt2\}$, $b_3\in\{2,2.51,2\sqrt2\}$,
$a_2,a_4\in\{2,2.51\}$.  Assume that the diagonal between corners
$1$ and $3$ has length in $[2\sqrt2,3.2]$, and that the other
diagonal has length $\ge3.2$.  Let $k_0$, $k_1$, $k_2$ be the
number of $b_i$ equal to $2$, $2.51$, $2\sqrt2$, respectively.
If $b_4=2.51$ and  $b_3=2$, no such subcluster exists
(the reader can check that $\Delta(4,4,x_3,4,2.51^2,x_6)<0$
under these conditions),
and we exclude this case.

1: $\vor_0(Q)< Z(4,k_1+k_2) - 0.009 k_2 -(k_0+2k_2)0.008/3$.
    \refno{193776341}

2: $\tau_0(Q) > D(4,k_1+k_2) + 0.04683 +
(k_0+2k_2-3)0.008/3+0.0066k_2$.
    \refno{898647773}
\smallskip

3: $\vor_0(Q) < Z(4,2) -0.0461 - 0.009 -2(0.008)$,
    if $a_2\in\{2,2.51\}$, $a_4=2$,
    $b_4=2\sqrt2$, $b_3=2.51$ or $2\sqrt2$.
    \refno{844634710}

4: $\tau_0(Q) > D(5,1) +0.04683+0.008+2(0.0066)$,
    if $a_2\in\{2,2.51\}$, $a_4=2$,
    $b_4=2\sqrt2$, $b_3=2.51$ or $2\sqrt2$.
    \refno{328845176}

5: $\vor_0(Q) < s_5 - 0.0461-0.008$,
    if $a_2\in\{2,2.51\}$, $a_4=2$,
    $b_3=2$,
    $b_4=2\sqrt2$.
    \refno{233273785}

6: $\tau_0(Q) > t_5 +0.008$,
    if $a_2\in\{2,2.51\}$,
    $a_4=2$,
    $b_3=2$,
    $b_4=2\sqrt2$.
    \refno{966955550}

(The penalties used in $\A_{20}$ are from Sections \refz 5.4 and
\refz  5.5.)

\subhead Section $\A_{21}$\dag\endsubhead
Recall that $\maxpi=0.06688$.

1: $\vor_0(S(2,2,2,y_4,2,2))+\vor_0(S(2,2,2,y_4',2,2))
    +\vor_0(S(2,2,2,y_4,y_4',2))
        < s_5 -0.008$, if $y_4,y_4'\in[2\sqrt2,3.2]$.
    \refno{275286804}

2: $\tau_0(S(2,2,2,y_4,2,2))+\tau_0(S(2,2,2,y_4',2,2))
    +\tau_0(S(2,2,2,y_4,y_4',2))
        > t_5+ 0.008$, if $y_4,y_4'\in[2\sqrt2,3.2]$.
    \refno{627654828}

3: $\vor_0(S(2,2,2,y_4,y_5,y_6)) < -2(0.008) + s_6-3(0.0461)$, if
    $y_4,y_5,y_6\in[2\sqrt2,3.2]$. (Compare $\A_{13}$.)
    \refno{995177961}

4: $\tau_0(S(2,2,2,y_4,y_5,y_6)) > t_6 + \maxpi$, if
    $y_4,y_5,y_6\in[2\sqrt2,3.2]$.
    \refno{735892048}

\subhead Section $\A_{22}$\dag\endsubhead  In $\A_{22}$ and $\A_{23}$,
$y_1\in [2.51,2\sqrt2]$, $y_4\in[2\sqrt2,3.2]$, and
$\dih<2.46$.
$\vor_0(Q)$ denotes the truncated Voronoi function
on the union of an anchored simplex and an adjacent special simplex.
Let $S'$ be the special simplex.  By deformations,
$y_1(S')\in\{2,2.51\}$.  If $y_1(S')=2.51$, the verifications
follow from $\A_6$ and $\vor_0(S')\le0$.  We may assume that
$y_1(S')=2$.  Also by deformations, $y_5(S')=y_6(S')=2$.

1:  $\vor_0(Q) < -3.58 + 2.28501 \dih$
    \refno{53502142}

2: $\vor_0(Q) < -2.715 + 1.67382 \dih$
    \refno{134398524}

3: $\vor_0(Q) < -1.517 + 0.8285 \dih$
    \refno{371491817}

4: $\vor_0(Q) < -0.858 + 0.390925 \dih$
    \refno{832922998}

5: $\vor_0(Q) < -0.358+0.009 + 0.12012 \dih$
    \refno{724796759}

6: $\vor_0(Q) < -0.186+0.009 + 0.0501 \dih$
    \refno{431940343}

When the cross-diagonal drops to $2.51$. We break $Q$ into two simplices
in the other direction.
Let $S''$ be an upright quarter with $y_5=2.51$.  In the next group
$\vor_0=\vor_0(S'')$

7:  $\vor_0 < -3.58/2 + 2.28501 \dih$
    \refno{980721294}

8: $\vor_0 < -2.715/2 + 1.67382 \dih$
    \refno{989564937}

9: $\vor_0 < -1.517/2 + 0.8285 \dih$
    \refno{263355808}

10: $\vor_0 < -0.858/2 + 0.390925 \dih$
    \refno{445132132}

11: $\vor_0 < (-0.358+0.009)/2 + 0.12012 \dih+0.2(\dih-1.23)$
    \refno{806767374}

12: $\vor_0 < (-0.186+0.009)/2 + 0.0501 \dih+0.2(\dih-1.23)$
    \refno{511038592}

\subhead Section $\A_{23}$\dag\endsubhead
$\tau_0(Q)$ denotes the truncated Voronoi function
on the union of an anchored simplex
(with $y_1\in[2.51,2\sqrt2]$, $y_4\in[2\sqrt2,3.2]$, $\dih<2.46$)
and an adjacent special simplex.

1: $-\tau_0(Q)+0.06585 < -3.48 + 2.1747 \dih$
    \refno{4591018}

2: $-\tau_0(Q)+0.06585  < -3.06 + 1.87427 \dih$
    \refno{193728878}

3: $-\tau_0(Q)+0.06585  < -1.58 + 0.83046 \dih$
    \refno{2724096}

4: $-\tau_0(Q)+0.06585  < -1.06 + 0.48263 \dih$
    \refno{213514168}

5: $-\tau_0(Q)+0.06585  < -0.83 + 0.34833 \dih$
    \refno{750768322}

6: $-\tau_0(Q)+0.06585  < -0.50 + 0.1694 \dih$
    \refno{371464244}

7: $-\tau_0(Q)+0.06585  < -0.29 +0.0014 + 0.0822 \dih$
    \refno{657011065}

Let $S'$ be the special simplex.  By deformations, we have
$y_5(S')=y_6(S')=2$, and $y_1(S')\in\{2,2.51\}$.  If
$y_1(S')=2.51$, and $y_4(S')\le3$, the inequalities listed
above follow from Section $\A_7$ and the inequality

8:  $\tau_0(S') > 0.06585$, if $y_1=2.51$, $y_4\in[2\sqrt2,3]$,
    $y_5=y_6=2$.
    \refno{66753311}

Similarly, the result follows if $y_2$ or $y_3\ge2.2$ from the
inequality

9:    $\tau_0(S') > 0.06585$, if $y_4\in[3,3.2]$, $y_5=y_6=2$,
        $y_1=2.51$, $y_2\in[2.2,2.51]$.
    \refno{762922223}

\medskip
Because of these reductions, we may assume in the first
batch of inequalities of $\A_{23}$ that when $y_1(S')\ne2$,
we have that $y_1(S')=2.51$, $y_5(S')=y_6(S')=2$, $y_4\in[3,3.2]$,
$y_2(S'),y_3(S')\le2.2$.  In all but {\tt (371464244)} and
{\tt (657011065)}, if $y_1(S')=2.51$, we prove the inequality
with $\tau_0(S')$ replaced with its lower bound $0$.

\smallskip
Again if the cross-diagonal is $2.51$, we break $Q$ in the other
direction.
Let $S''$ be an upright quarter with $y_5=2.51$.
Set $\tau_0 = \tau_0(S'')$. We have

10: $-\tau_0+0.06585/2 < -3.48/2 + 2.1747 \dih$
    \refno{953023504}

11: $-\tau_0+0.06585/2  < -3.06/2 + 1.87427 \dih$
    \refno{887276655}

12: $-\tau_0+0.06585/2  < -1.58/2 + 0.83046 \dih$
    \refno{246315515}

13: $-\tau_0+0.06585/2  < -1.06/2 + 0.48263 \dih$
    \refno{784421604}

14: $-\tau_0+0.06585/2  < -0.83/2 + 0.34833 \dih$
    \refno{258632246}

15: $-\tau_0+0.06585/2  < -0.50/2 + 0.1694 \dih+0.03(\dih-1.23)$
    \refno{404164527}

16: $-\tau_0+0.06585/2  < -0.29/2 +0.0014/2 + 0.0822 \dih
+0.2(\dih-1.23)$
    \refno{163088471}

\bigskip
\subhead Section $\A_{24}$\dag\endsubhead
These final calculations here are used to determine what is squandered
when $\dih>2.46$.

1: $\tau_0+0.0822\dih> 0.159$, if $y_1\in[2.51,2\sqrt2]$,
$y_6\in[2.51,2.75]$, $y_2=y_4=2$. \refno{968721007}

2: $\dih<1.23$ if $y_1\in[2.51,2\sqrt2]$, $y_6\ge2.51$,
$y_2=2.51$, $y_4=2$. \refno{783968228}

3: $\dih<1.23$ if $y_1\in[2.51,2\sqrt2]$, $y_6\ge2.75$,
$y_2=y_4=2$. \refno{745174731}

\heads{Appendix 2}
\vfill\eject
\head Appendix 2. Some Conventions\endhead

\parindent=10pt

Throughout the paper, we have preferred to work with
compact domains.  As we divide cases into
compact sets, boundaries will overlap.
This leads to various mild inconsistencies unless certain
statements in the paper are interpreted appropriately.

For example, if an edge of a quasi-regular tetrahedron
is exactly 2.51, the quasi-regular tetrahedron is also
a quarter.  If some of the simplices along that edge are
interpreted as quasi-regular tetrahedra and others are
interpreted as quarters, this could easily have unintended
effects.  In such cases we ask the reader to decide once
and for all whether the edge is to be considered the diagonal
of a quarter or as a short edge of a quasi-regular tetrahedron,
and then adhere to that convention.

In general when a length lies on the boundary between
two cases,  the inequalities have been designed
to hold for whichever of the two cases is selected, as long as
the selection is consistently adhered to.

When we divide the domain into several compact regions, and divide
a function piecewise on each region,  in several places we use an
abbreviated style that might create ambiguities for function
values at boundary cases.  Again we ask the reader to adhere to
any consistent convention.

\smallskip
In most cases, bounds on the score are strict.
There are only a few places where exact equality can be obtained
and where it makes an appreciable difference.
The most significant are the bounds $\sigma\le\pt$ on
quasi-regular tetrahedra and $\sigma\le0$ on quad-clusters.
The fact that these are attained for the regular cases with edge
lengths 2 and diagonal $2\sqrt{2}$ on the quad-cluster and for
no other cases gives the bound $\pi/\sqrt{18}$ on density and
the local optimality of the fcc and hcp packings.

Another place where we have allowed equality to be obtained is
with $\tau_0\ge0$ for quasi-regular simplices.  The importance
of equality for Rogers's bound on the density of packings is explained
in III.

There are also a few less significant cases where an inequality
is sharp.
For example, $$\tau_0(2.51,2,2,x,2,2)\ge0, \quad\vor_0(2.51,2,2,x,2,2)\le0$$
for special simplices satisfying
 $x\in[2\sqrt{2},3.2]$.  Also, equality occurs in Lemma
F.1.9 and F.2.2.

\bye